
\input amstex
\documentstyle{amsppt}
\magnification=1200

\redefine\R{{{\bold R}}}
\redefine\Z{{{\bold Z}}}
\redefine\Q{{{\bold Q}}}
\redefine \Cplx{{{\bold C}}}

\define\ab{\allowbreak}

\define\hy{{\hbox{\it -}}}
\define\one{{1\hy}}
\define\hyspec{{{\hbox{\it -spec}}}}
\define\spec{\text{\it spec}}
\define\Vol{\text{\it Vol }}

\define\mm{{{\frak g}}}
\define\hh{{{\frak h}}}
\define\nn{{{\frak N}}}
\define\cc{{{\frak C}}}
\redefine\aa{{{\frak A}}}
\redefine\gg{{{\frak g}}}
\define\zz{{{\frak z}}}

\define\Remark{{{\flushpar {\bf Remark.}\quad}}}

\define\F{{F_\tau}}

\define\p#1#2{{\frac{\partial}{\partial {#1}_{#2}}}}
\define\ps#1#2{{\frac{\partial^2}{\partial {#1}^2_{#2}}}}

\define\m#1{{{\gg}^{(#1)}}}
\define\h#1{{{\hh}^{(#1)}}}
\define\g#1{{{\gg}^{(#1)}}}

\define\BX{{\bar{X}}}
\define\BY{{\bar{Y}}}
\define\BZ{{\bar{Z}}}
\define\BE{{\bar{E}}}


\define\Tau{{{\Cal T}}}

\redefine\L{{{\lambda}}}

\redefine\H{{{\Cal H}}}
\redefine\xi{{{\zz}}}

\define\Gup#1{{G^{(#1)}}}
\define\GGup#1{{\GG^{(#1)}}}
\define\GGperpup#1{{\GG_\perp^{(#1)}}}
\define\Gbar{{\bar{G}}}

\define\M#1{{G^{(#1)}}}
\define\gggbar{{\bar{\gg}}}
\define\GG{{\Gamma}}
\define\G#1{{{{\Gamma}_{#1}}}}
\define\GGbar{{\bar{\Gamma}}}
\define\GGGbar#1{{{{\bar{\Gamma}}_{#1}}}}

\define\mbar{\gggbar}
\define\Mbar{\Gbar}

\define\Rep{{\raise 1.5pt\hbox{$\rho$}_\GG}}
\define\rep#1{{\raise 1.5pt\hbox{$\rho$}_\G{#1}}}

\define\Nilmfld{{\GG \backslash G}}
\define\nilmfld#1{{\G{#1} \backslash G}}
\define\Nilmfldbar{{\GGbar \backslash \bar{G}}}
\define\nilmfldbar#1{{\GGGbar{#1} \backslash \bar{G}}}

\define\gbar{{\bar{g}}}

\define\Fnsp{{L^2(\Nilmfld)}}
\define\fnsp#1{{L^2(\nilmfld{#1})}}
\define\fnspbar#1{{L^2(\nilmfldbar{#1})}}
\define\Fnspg{{L^2(\GG \backslash G)}}
\define\fnspg#1{{L^2(\G{#1} \backslash G)}}

\redefine\d#1#2{{\dot{#1}_{#2}(s)}}\redefine\dL#1#2{{\dot{#1}_{#2}(s+\L)}}

\define\dsz#1#2{{{\frac{d{#1}_{#2}}{ds}}_{|_0}}}
\redefine\ss#1#2{{<\sss,{#1}_{#2}>}}
\redefine\b#1#2{{\overline{{#1}}_{#2}}}

\define\endpf{\hbox{\vrule height1.5ex width.5em}}


\topmatter

\pagewidth{5.5 in}
\pageheight{7.5 in}
\hoffset=-.02 in

\title 
A New Construction of Isospectral Riemannian Nilmanifolds with Examples
\endtitle

\rightheadtext{A New Construction of Isospectral Nilmanifolds}

\author
Ruth Gornet
\endauthor

\address
Ruth Gornet; Texas Tech University; Department of Mathematics; 
Lubbock, Texas \ 79409-1042; email: gornet\@math.ttu.edu
\indent September 1994
\endaddress

\keywords
Isospectral manifolds; Spectrum of the Laplacian on forms; 
Nilpotent Lie groups;  Cocompact, discrete subgroups; 
Quasi-regular representations; Square-in\-te\-gra\-ble 
representations.
\flushpar\indent Research at MSRI supported in part by NSF
grant DMS-9022140. Research at MSRI and Texas Tech supported in part
by NSF grant DMS-9409209
\endkeywords


\abstract
We present a new construction for obtaining pairs of higher-step 
isospectral Riemannian nilmanifolds and compare several resulting new examples.
In particular, we present new examples of manifolds that are isospectral
on functions, but not isospectral on one-forms. 
\endabstract
\endtopmatter

\document
\subheading{\S1 Introduction}
\bigskip

The {\it spectrum} of a closed Riemannian manifold $(M,g)$, 
denoted $\spec(M,g)$, is the collection of eigenvalues with multiplicities of the 
associated La\-place--Bel\-tra\-mi operator acting on smooth functions.
Two Riemannian manifolds $(M,g)$ and $(M',g')$ are said to be
{\it isospectral} if $\spec(M,g)=\spec(M',g').$
A basic question in spectral geometry is determining what geometric 
information is contained in the spectrum of a Riemannian manifold.
\medskip

Despite considerable research in the area, only a few geometric 
properties are known to be spectrally determined; for example, 
dimension, volume, and total scalar curvature.
Examples of isospectral manifolds provide us with the only means for determining 
properties not determined by the spectrum.  
\medskip
 
Most previous constructions for producing isospectral manifolds are 
based solely on representation theory.  In fact, most known examples 
can be explained by Sunada's method \cite{S} or its generalizations 
\cite{GW1}, \ab \cite{DG2}. (See \cite{B2} for a general overview.)
All manifolds constructed by these methods are
{\it strongly isospectral;} that is, all natural, strongly elliptic, 
self-adjoint operators on the manifolds are isospectral.  In particular, 
these manifolds share the same $p$-form spectrum.
The {\it $p$-form spectrum} of a manifold is the eigenvalue spectrum of the 
Laplace-Beltrami operator extended to act on smooth $p$-forms for $p$ a positive 
integer.
\medskip

The primary goal of this paper is the development of a new construction 
for producing pairs of isospectral nilmanifolds. Two-step nilmanifolds have been a 
rich source of examples of isospectral manifolds.  
Moreover, their  geometry has been studied in some detail.
(See \cite{DG1}, \ab \cite{BG}, \ab \cite{O}, \ab \cite{P1}, \ab \cite{P2}, \ab \cite{P3},
\ab \cite{GW1}, \ab \cite{GW2}, \ab \cite {G1}, \ab \cite {G2}, \ab and  \cite {E}.)  
This new construction uses 
techniques from Riemannian geometry, Lie groups, and representation theory 
to produce pairs of isospectral nilmanifolds of arbitrary step.  The higher-step 
examples have a much richer geometry, exhibiting many properties 
not previously found.
\medskip

While representation theory is used as a tool, this construction differs 
from previous ones in that the resulting pairs of isospectral manifolds 
need not be isospectral on one-forms, and so do not fall under the traditional
Sunada setup.   This property was previously exhibited by pairs 
of isospectral Heisenberg manifolds constructed by Gordon and Wilson \cite{GW2},
\cite{G2}.  And for any choice of $P$, Ikeda \cite{I2} has constructed 
examples of isospectral lens spaces that are isospectral on $p$-forms for 
$p=0,1,\cdots,P$ but not isospectral on $(P+1)$-forms. These are the only 
known examples.  
\medskip

The only other examples of isospectral manifolds that do not fall under the 
traditional Sunada construction are bounded domains of Urakwa \cite{U},  and nonlocally 
isometric examples of Szabo \cite{Sz2}, Gordon \cite{G3}, \cite{G4}, and 
Gordon--Wilson \cite{GW3}.
\medskip

Recent results of Pesce \cite{P4} now explain the Ikeda and Urakawa examples in 
a Sunada-like setting. This setting requires a genericity  assumption that 
excludes nilmanifolds. Moreover, the construction presented here 
generalizes the method used by Gordon--Wilson to construct the Heisenberg examples.
\medskip

Consequently, outside of the nonlocally isometric examples mentioned above,  
the construction below subsumes all known examples of isospectral manifolds 
that do not fall under a Sunada setup.
\medskip

After establishing notation in Section 2, the new construction for 
producing pairs of higher-step isospectral nilmanifolds is presented in Section 3. 
In Section 4, we use the new construction to 
produce 
pairs of isospectral three-step nilmanifolds with the following combinations of properties.
\medskip

\centerline{Table I:  New Examples of Isospectral Manifolds}
\medskip
{\eightpoint
\centerline{$
\vbox{\offinterlineskip
\halign{\strut\vrule#&\ #\hfil\ &&\vrule#&\ \hfil#\hfil\ \cr
\noalign{\hrule}
&Pair of 3-Step&&$\forall p$ Same&&Rep. Equiv.&&Isomorphic&&Same&&Same&\cr
&Isospectral&&$p$-form&&Fundamental&&Fundamental&&Length&&Marked Length&\cr
&Nilmanifolds&&Spectrum&&Groups&&Groups&&Spectrum&&Spectrum&\cr
\noalign{\hrule}
&I(7 dim)&&Yes&&Yes&&No&&No&&No&\cr
\noalign{\hrule}
&II(5 dim)&&Yes&&Yes&&Yes&&Yes&&No&\cr
\noalign{\hrule}
&III$\backslash$IV(7$\backslash$5 dim)&&No&&No&&No&&No&&No&\cr
\noalign{\hrule}
&V(7 dim) &&No&&No&&Yes&&Yes&&Yes&\cr
\noalign{\hrule}}}
$}
}
\bigskip

The properties listed in the above table are defined as follows.
Two cocompact (i.e. $\Nilmfld$ compact), discrete subgroups $\Gamma_1$ and 
$\Gamma_2$ of a Lie group $G$ are called {\it representation equivalent} if 
the associated quasi-regular representations are unitarily equivalent. 
(See Section 2 for details.)
The {\it length spectrum} of a Riemannian manifold is the set of 
lengths of closed geodesics, counted with multiplicity.  The multiplicity of
a length is defined as the number of distinct free
homotopy classes in which the length occurs.  (Note:  other definitions of 
multiplicity also appear in the literature.)  The pairs of isospectral 
manifolds above have the same lengths of closed geodesics. 
However, the length spectra often differ in the multiplicities that occur.
The {\it marked length spectrum} takes into account not only the lengths of 
the closed geodesics but also the free homotopy classes in which the geodesics
occur.
\medskip

In Section 4 we compare the quasi-regular representations 
and the fundamental groups of Examples I through V. 
We also compare the $p$-form spectrum, but the calculations are left
to an Appendix.  The length spectrum and marked 
length spectrum of these examples will be examined in \cite{Gt4}.
\medskip

Example I is the first example of a pair of {\it nonisomorphic,}
representation equivalent,  cocompact, discrete subgroups of a nilpotent 
Lie group.  It is also the first example of a pair of representation 
equivalent cocompact, discrete subgroups of a solvable Lie group
producing Riemannian manifolds 
that do not have the same length spectrum.
This example has implications in representation theory on nilpotent Lie groups  
and motivated  \cite{Gt1} and \cite{Gt2}.
\medskip

We prove in the Appendix that the manifolds in Examples III, IV and V are not 
isospectral on one-forms.
Outside of the traditional Sunada setup, no general method is known for comparing the 
one-form spectrum of manifolds.
The methods illustrated in the Appendix are new, 
as previously used
techniques could not be applied to the higher-step examples. 
The only previous examples of manifolds that are isospectral on functions but 
not isospectral on $p$-forms for all $p$
are the lens spaces and Heisenberg manifolds mentioned above.
\medskip

Example V is the first example of a pair of Riemannian manifolds with the 
same marked length spectrum, but not the same spectrum on one-forms.  This 
example contrasts with two-step results relating the marked length 
spectrum and the $p$-form spectrum \cite{E}.  This example will be studied
in detail in \cite{Gt4}.
\medskip

Most of the contents of this paper are contained in the author's thesis
at Washington University in St. Louis in partial fulfillment of the requirements
for the degree of Doctor of Philosophy.  The author wishes to express
deep gratitude to her advisor, Carolyn S. Gordon, for all of her suggestions,
encouragement and support.
\bigskip

\bigskip


\bigskip
\subheading{\S2 Background and Notation}
\bigskip

Let $G$ be a simply connected Lie group and let $\Gamma$ be a cocompact, 
discrete subgroup of $G.$ A Riemannian metric $g$ is {\it left invariant} 
if the left translations of $G$ are isometries.   The left invariant metric 
$g$ projects to a Riemannian metric on $\Nilmfld,$ which we also denote by 
$g.$  Note that a left invariant metric is determined by a choice of 
orthonormal basis of the Lie algebra $\gg$ of $G.$
\medskip

As $G$ is unimodular, the Laplace--Beltrami operator of $(\Nilmfld, g)$ may be 
written 
$$\Delta = -\sum_{i=1}^n {E_i}^2,\tag2.1$$ 
where $\{E_1, \cdots, E_n\}$ is an orthonormal basis of the Lie algebra $\gg$ 
of $G.$
\medskip 

The Laplace-Beltrami operator acting on smooth $p$-forms is defined by 
$\Delta = d \delta + \delta d.$  Here $\delta$ is the metric adjoint of $d.$  
Equivalently $\delta = (-1)^{n(p+1)+1}*d*,$ where $*$ is the Hodge-$*$ operator.  
We denote the {\it $p$-form spectrum} of a Riemannian manifold $(M,g)$  by 
$p\hy\spec(M,g).$
\medskip

The quasi-regular representation $\Rep$ of $G$ on $\Fnspg$ is 
defined as follows: 

\flushpar for all $x$ in $G$ and  $f$ in $\Fnspg,$ 
$$\Rep(x)f=f\circ R_x.$$
Here $R_x$ denotes the right action of $x$ on $\GG \backslash G.$
The quasi-regular representation is known to be unitary.
\medskip

We say $\G1$ and $\G2$ are {\it representation equivalent} if $\rep1$ 
and $\rep2$ are unitarily equivalent; that is, $\G1$ and $\G2$ are 
representation equivalent if there exists a unitary isomorphism 
$T:\fnspg1 \rightarrow \fnspg2$  such that $T(\rep1(x)f)= \rep2(x)Tf$ 
for every $x$ in $G$ and every $f$ in $\fnspg1.$
\bigskip

\proclaim{Proposition 2.2 (Gordon--Wilson \cite{GW1}) }
Let $\G1$ and $\G2$ be cocompact, discrete subgroups of a simply connected 
Lie group $G.$  Let $g$ be a left invariant metric on $G.$  
If $\G1$ and $\G2$ are representation equivalent, then 
$$p\hy \spec(\G1\backslash G,g) = p\hy \spec(\G2 \backslash G,g)$$
for $p=0,1, \cdots, dim(G).$ 
\endproclaim
\bigskip

\Remark Pairs of isospectral manifolds constructed using the traditional Su\-na\-da 
method are of the form $(\Gamma_1 \backslash M,g)$ and 
$(\Gamma_2 \backslash M,g)$ 
where $\Gamma_1$ and $\Gamma_2$ are representation equivalent, cocompact, 
discrete subgroups of a group $G$ acting by isometries on 
a Riemannian manifold $(M,g).$
\bigskip

For a Lie algebra $\gg$, denote by $\g1$  
the derived algebra $[\gg,\gg]$ of $\gg.$ 
That is, $\g1$ is the Lie subalgebra of $\mm$ generated
by all elements of the form $[X,Y]$ for $X,Y$ in $\gg.$
Inductively, define $\g{k+1}=[\gg,\g{k}].$
A Lie algebra
$\mm$ is said to be {\it k-step nilpotent} if $\m{k} \equiv 0$
but $\m{k-1} \not \equiv 0.$
A Lie group $G$ is called {\it $k$-step nilpotent} if its Lie algebra is.
\medskip

Let $\Gup{k} = exp(\m{k})$ denote the $k$th derived subgroup of $G.$ 
We denote the center of $G$ by $Z(G)$ and the center of $\mm$ by $\zz.$
Note that if $G$ is $k$-step nilpotent, then $\Gup{k-1} \subset Z(G).$
\medskip

Let $exp$  denote  the Lie algebra 
exponential from $\mm$ to $G.$
The Camp\-bell--Baker--Haus\-dorff formula gives us the group operation of 
$G$ in terms of $\mm.$  Namely, for $X,Y \in \mm:$
$$exp(X)exp(Y) = exp(X+Y+\frac1{2}[X,Y]+\frac1{12}[X,[X,Y]]+
\frac1{12}[Y,[Y,X]]+\cdots),$$
where the remaining terms are higher-order brackets.
Note that for two-step nilpotent Lie groups, only the first three terms in the
right-hand side are nonzero.  For three-step groups, only the first
five terms are nonzero.
If $\mm$  is nilpotent and $G$ is simply connected,
then $exp$ is a diffeomorphism from $\mm$
onto $G.$  Denote its inverse by $log.$
\medskip

Let $\G1$ and $\G2$ be cocompact, discrete subgroups of nilpotent Lie 
groups $G_1$ and $G_2$ respectively.  Any abstract group isomorphism 
$\Phi:\G1 \rightarrow \G2$ extends uniquely to a Lie group isomorphism 
$\Phi:G_1 \rightarrow G_2.$
\medskip

Let $\GG$ be a cocompact,
discrete subgroup of a nilpotent Lie group $G$ with left invariant metric
$g.$   
The locally homogeneous space $(\Nilmfld,g)$ is called a {\it 
Riemannian nilmanifold.}
If $G$ is an abelian Lie group,
then $\GG$ is merely a lattice of full rank in $G,$
and in this case $log\GG$ is also a lattice in $\mm.$
\medskip

Let $\mm_\Q= span_\Q \{ log \GG \}.$
This is a {\it rational} Lie algebra; that is, 
there exists a basis of $\mm$ made up of elements of $log\GG$
such that the structure constants are rational.
A Lie subalgebra  $\hh$ of $\mm$ is called
{\it rational} if $\hh$ is spanned by
$\hh \cap \mm_{\bold Q}.$
Note that the notion of rational depends on $\GG.$
If $H=exp(\hh)$ is the connected Lie subgroup of $G$ 
with rational Lie algebra $\hh,$
then $\GG \cap H$ is a cocompact, discrete
subgroup of $H.$
The $\m{k}$ are always rational Lie subalgebras of $\mm.$
\medskip

The Kirillov theory of  irreducible unitary representations
of nilpotent groups gives us a correspondence between 
irreducible unitary representations of $G$
and elements of  $\mm^*,$
the dual of $\mm.$
In particular, fix $\tau \in \mm^*.$
Let $\hh$ be a rational subalgebra of $\mm$
that is maximal with respect to the property that
$\tau([\hh,\hh])\equiv 0.$
The subalgebra $\hh$ is called a {\it polarization} of $\tau.$
Let $H=exp(\hh)$ be the connected subgroup of
$G$ with Lie algebra $\hh.$
\medskip

Define a character $\overline{\tau}$ of $H$ by 
$${\overline{\tau}}(h)=e^{2\pi i \tau (log(h))}\tag2.3$$
for all $h$ in $H.$
Define $\pi_\tau$ to be the irreducible representation
of $G$ induced by the representation $\overline{\tau}$ of $H.$
Denote by  ${\H}_\tau$ the representation space of $\pi_\tau.$
Two such irreducible representations
$\pi_\tau$ and $\pi_{\tau'}$ are unitarily
equivalent if and only if $\tau' = \tau \circ Ad(x)$
for some $x$ in $G.$  Here $Ad(x)$ 
is the adjoint map from $\mm$ to $\mm.$
\medskip

For $\tau$ in $\mm^*,$  the {\it coadjoint orbit\/} of $\tau$ is
$$O(\tau) = \{ \tau \circ Ad(x): x \in G \}.$$
Hence $\pi_\tau$ and $\pi_{\tau'}$ are unitarily
equivalent if and only if $\tau$ and $\tau'$ lie in the same coadjoint 
orbit of $\mm^*.$
\medskip

As $G$ is nilpotent, every irreducible representation of $G$
is unitarily equivalent to $\pi_\tau$
for some $\tau \in \mm^*,$ and the quasi-regular representation $\Rep$ is 
completely reducible.
Thus the representation space $\Fnsp$ is unitarily isomorphic to 
$$\Fnsp \cong {\underset{\tau \in \Tau}\to{\bigoplus}}m(\tau){\H}_\tau$$
for some $\Tau \subset \mm^*.$  
Here $m(\tau)$ denotes the multiplicity of $H_\tau,$ and we assume $\Tau$ 
contains at most one element of each coadjoint orbit of $\mm^*.$
\medskip

A good reference for representation theory on nilpotent Lie groups is
\cite{CG}.

\bigskip

\bigskip

\subheading{\S3 A New Construction of Isospectral Nilmanifolds}

\bigskip

Let $G$ be a simply connected, $k$-step nilpotent Lie group with Lie algebra $\mm.$  
Define $\Gbar$ to be the simply connected, $(k$--$1)$-step nilpotent Lie group 
$G/\Gup{k-1}.$  For $\Gamma$ a cocompact, discrete subgroup of $G,$
denote by $\bar\Gamma$ the image of $\Gamma$ under the canonical projection 
from $G$ onto $\Gbar.$
The group $\bar\Gamma$ is then a cocompact, discrete subgroup of $\Gbar.$
For a left invariant metric $g$ on $G,$ we associate a left invariant 
metric $\gbar$ on $\Gbar$ by restricting $g$ to an orthogonal complement 
of $\g{k-1}$ in $\mm.$
\medskip

We call the $(k$--$1)$-step nilmanifold $(\Nilmfldbar,\gbar)$ the {\it quotient 
nilmanifold of 
$(\Nilmfld,g).$}  By using the definition of $\gbar,$ one easily sees that the projection 
$(\Nilmfld, g) \rightarrow (\Nilmfldbar,\gbar)$ is a Riemannian submersion.
\medskip

The Lie algebra $\gggbar$ of $\Gbar$ is just $\mm / \g{k-1}.$   We denote elements 
of $\gggbar$ by $\bar{U},$ where $\bar{U}$ is the image of $U$ under the canonical 
projection from $\mm$ onto $\gggbar.$
\bigskip

\proclaim{Definition 3.1} Let $G$ be a simply connected nilpotent Lie group.  We 
say $G$ is {\it strictly nonsingular} if the following property holds:  for 
every $z$ in $Z(G)$ and every noncentral $x$ in $G$ there exists an element $a$ 
in $G$ such that $$[x,a] = z.$$
Here $[x,a]$ denotes the commutator $[x,a]=xax^{-1}a^{-1}.$

Equivalently, the Lie algebra $\gg$ is {\it strictly nonsingular} if for every 
noncentral $X$ in $\gg,$  
$$\zz \subset ad(X)(\gg).$$
That is, for every $X$ in $\gg-\zz$ and every $Z$ in $\zz$ there exists a 
vector $Y$ in $\gg$ such that $[X,Y]=Z.$
\endproclaim
\bigskip

\proclaim{Theorem 3.2}
Let $G$ be a simply connected, strictly nonsingular nilpotent Lie group with left 
invariant metric $g.$  If \ $\G1$ and $\G2$ are cocompact, discrete subgroups 
of $G$ such that 
$$\G1 \cap Z(G) = \G2 \cap Z(G) \quad \text{and} \quad \spec(\nilmfldbar1,\gbar) = \spec(\nilmfldbar2,\gbar),$$ then  $$\spec(\nilmfld1,g) = \spec(\nilmfld2,g).$$
\endproclaim
\bigskip

\Remark The above construction is a generalization of the construction used by Gordon 
and Wilson to obtain pairs of isospectral Heisenberg manifolds  
\cite {GW2}.  If we let the Lie group $G$ be a simply connected, strictly 
nonsingular, two-step nilpotent Lie group with a one-dimensional center, 
then $G=H_n$ for some $n,$ where $H_n$ denotes the $(2n+1)$-dimensional Heisenberg group.  
\bigskip
\pagebreak

\demo{Proof of Theorem 3.2}
\bigskip

We use the notation of Section 2.
\medskip

For $i=1,2,$ let $\Tau_i$  be a subset of $\gg^*$  such that 
$$\fnsp{i} \cong {\underset{\tau\in\Tau_i}\to\bigoplus} {m_i(\tau)}\H_\tau.$$
Recall that $m_i(\tau)$ denotes the multiplicity of $\pi_\tau$ in the quasi-regular 
representation of $G$ on $\fnsp{i},$ and we assume that $\Tau_i$ contains at 
most one element of each coadjoint orbit of $\gg^*.$
\medskip

We decompose the index set $\Tau_i =\Tau'_i \cup \Tau''_i$ by letting 
$$\Tau'_i = \{\tau \in \Tau_i :  \tau(\zz) \equiv 0\}, \quad \text{and} \quad  
\Tau''_i = \{\tau \in \Tau_i :  \tau(\zz )\not \equiv 0\}.$$
We likewise decompose the representation space $\fnsp{i}$ by letting  
$$\H'_i= {\underset{\tau\in\Tau'_i}\to\bigoplus} m_i(\tau)\H_\tau \quad \text{and} \quad \H''_i= {\underset{\tau \in \Tau''_i}\to\bigoplus} m_i(\tau)\H_\tau .$$
As representation spaces 
$\fnsp{i}= \H'_i \oplus \H''_i.$
\medskip

Clearly
$spec(\nilmfld{i}, g) = spec'(\nilmfld{i}, g) \cup spec''(\nilmfld{i}, g),$ 
where $spec'(\nilmfld{i}, g)$ and 
\break
$spec''(\nilmfld{i}, g)$ are defined as the 
spectrum of the Laplacian restricted to acting on $\H'_i,$ and  $\H''_i,$  
respectively.   The multiplicity of an eigenvalue in 
$spec(\nilmfld{i}, g)$ is equal to the sum of its multiplicities 
in  $spec'(\nilmfld{i}, g)$ and $spec''(\nilmfld{i}, g).$

\bigskip

\proclaim{Lemma 3.3}
The Laplacian of $(\nilmfld{i}, g)$ acting on $\H'_i$ is precisely the Laplacian of $(\nilmfldbar{i},\gbar)$ acting on $\fnspbar{i}.$
Thus $\spec(\nilmfldbar{i},\gbar) = \spec'(\nilmfld{i},g)$ for $i=1,2.$ 
\endproclaim
\bigskip

\proclaim{Lemma 3.4}
The representations of $G$
on $\H''_1$ and $\H''_2$ are unitarily equivalent, 
hence $\spec''(\nilmfld1,g) = \spec''(\nilmfld2,g).$
\endproclaim
\bigskip

Theorem 3.2 now follows.
\medskip

The proof of Lemma 3.3 is essentially an extension of the first part of 
the proof used by Gordon--Wilson to construct pairs of isospectral 
Heisenberg manifolds.  (See \cite{GW2}, Theorem 4.1)  The details are 
included here for completeness and because of a difference in notation.
\bigskip
\pagebreak

\demo{Proof of Lemma 3.3}
\bigskip

Let $\{Z_1, Z_2, \cdots, Z_T\}$ be an orthonormal basis of $\zz = \g{k-1}.$  Extend it to 
$\{E_1, E_2, \allowbreak \cdots, E_N,\allowbreak Z_1, Z_2, \cdots, Z_T\},$ an orthonormal basis of $\gg.$
By (2.1), the La\-place-Bel\-tra\-mi operator of $(G,g)$ is 
$$\Delta = - \sum_{n=1}^N{E_n}^2  - \sum_{k=1}^K{Z_k}^2.$$  
\medskip

View functions in $\fnsp{i}$ as left $\G{i}$-invariant functions of $G.$  The 
subspace $\H'_i$ is then those functions in $\fnsp{i}$ that are independent 
of the center, which correspond to functions in $\fnspbar{i}$ in a natural way.
So when we restrict $\Delta$ to $\H'_i,$ we have  
$$\Delta = - \sum_{n=1}^N{E_n}^2,$$ 
which corresponds to the Laplacian of $(\nilmfldbar{i},\gbar).$
\medskip

The Laplacian of $(\nilmfld{i}, g)$ acting on $\H'_i$ is then precisely the 
Laplacian of $(\nilmfldbar{i},\gbar)$ acting on $\fnspbar{i},$ 
so $spec(\nilmfldbar{i},\gbar) = spec'(\nilmfld{i},g),$
as 
desired.  $\endpf$
\enddemo

\bigskip

Before proving Lemma 3.4, we must introduce some of the theory of square 
integrable representations of nilpotent Lie groups.
\bigskip

\proclaim{Definition 3.5}
Let $G$ be a locally compact, unimodular group with center Z(G).  We say that an 
irreducible unitary representation $\pi$ of $G$ on a Hilbert space $\H$ 
is {\it square integrable\/} if there are nonzero vectors $x_1$ and $x_2$ 
in $\H$ such that
$$\int_{G/Z(G)} |(\pi(s)x_1,x_2)|^2 d\overline{\mu}(\overline{s}) <  \infty.$$
\endproclaim
\bigskip

Here $d\overline{\mu}(\overline{s})$ denotes integration over $G/Z(G)$ with 
respect to a choice of Haar measure $\overline{\mu}$ on $G/Z(G).$  As the 
center acts trivially, the integrand may be viewed as a function of $G/Z(G).$
\medskip

Let $\xi^\perp$ be the subalgebra of $\mm^*$ defined by 
$\xi^\perp = \{ \mu \in \mm^* : \mu(\xi) \equiv 0 \}.$  
Note that $\xi^\perp \cong (\mm/\xi)^*.$ For $\tau \in \mm^*,$ 
let $b_\tau$ denote the skew-symmetric, bilinear form on $\mm/\xi$ 
defined by $b_\tau(\BX,\BY)=\tau([X,Y])$ for all $\BX,\BY$ in $\mm/\xi.$  
Here $X$ and $Y$ are any elements of $\mm$ that project onto $\BX$ 
and $\BY$ respectively. 
\bigskip

\pagebreak

\proclaim{Theorem 3.6 (Moore--Wolf \cite {MW})}
For a linear functional $\tau$ in $\mm^*$ with coadjoint orbit $O(\tau)$ and 
corresponding irreducible unitary representation $\pi_\tau,$ the following 
three conditions are equivalent:

(1) $\pi_\tau$ is square integrable.

(2) $O(\tau) = \tau + \xi^\perp.$

(3) $b_\tau$ is nondegenerate on $\mm/\xi.$
\endproclaim
\bigskip


\demo{Proof of Lemma 3.4}
\bigskip

Fix  $\tau \in \Tau''_{i}.$ Let $Z \in \xi$ be such that $\tau(Z) \neq 0.$  
By strict nonsingularity, for all noncentral $X \in \mm$ there 
exists $Y \in \mm$ such that $[X,Y]=Z.$  Hence 
$b_\tau(\BX,\BY) = \tau([X,Y]) \neq 0.$  So $b_\tau$ is nondegenerate.
By Theorem 3.6, $\pi_\tau$ is square integrable.
Note that $b_\tau$ nondegenerate implies that $N=dim(\mm/\xi)$ is even.
\medskip

Recall from Section 2 that $\pi_\tau$ is independent of the choice of 
$\tau$ in $O(\tau).$ 
By Theorem 3.6, as $\pi_\tau$ is square integrable, 
the coadjoint orbit $O(\tau)$ is uniquely determined by the restriction of 
$\tau$ to the center.
We may thus assume  $\tau \in \xi^*.$
\medskip  

Let $\alpha$ be  a volume form  on $G/Z(G).$  That is, let $\alpha$ be a fixed, 
alternating, $N$-linear form over $\mbar = \mm/\xi.$
As $b_\tau$ is nondegenerate,
${b_\tau}^{(\frac1{2}N)}=b_\tau \wedge \cdots \wedge b_\tau$ is also a volume 
form on $G/Z(G)$ and hence a scalar multiple of $\alpha.$
Define $P_\alpha(\tau)$ by $$P_\alpha(\tau)\alpha = b_\tau \wedge \cdots \wedge b_\tau.$$
The polynomial $P_\alpha$ is homogeneous of degree $\frac1{2}N$ on $\mm^*$ and 
depends only on the choice of volume form $\alpha.$
\medskip

Let $\xi^*$ denote the dual of $\xi.$
Let $L=\G1 \cap Z(G) = \G2\cap Z(G)$
and let ${L}^*\subset \xi^*$ be the dual lattice of $L.$
\medskip

We now use the following occurrence and multiplicity condition, also due to Moore and Wolf.
\medskip

\proclaim{Theorem 3.7 (Moore--Wolf \cite {MW})}
Let $G$ be a nilpotent Lie group and $\GG$ a cocompact, discrete 
subgroup of $G.$  Let $L=\GG \cap Z(G).$  Fix a volume form 
$\alpha_\GG$ on $\mm/\xi$ so that $\GGbar \backslash \Mbar$ has 
volume one.  Let $\tau$ be a nonzero element of $\xi^*$ such 
that $\pi_\tau$ is square integrable.   The representation 
$\pi_\tau$ occurs in the quasi-regular representation of $G$ 
on $\Fnsp$ if and only if $\tau \in L^*.$  Moreover, its 
multiplicity $m(\tau)$ is $|P_{\alpha_\GG}(\tau)|.$
\endproclaim
\bigskip

By Theorem 3.7, the square integrable representation  $\pi_\tau$ 
occurs in the quasi-regular representation of $G$ on 
$\fnsp{i}$ if and only if $\tau$ is contained in ${L}^*.$  Thus 
for $i=1,2$ the coadjoint orbits represented in $\Tau''_i$  
correspond to the elements of ${L}^*.$   
We may assume $\Tau''_1 = \Tau''_2.$
\medskip

As the Riemannian metrics of $(\nilmfldbar1,\gbar)$ and 
$(\nilmfldbar2,\gbar)$ arise from the same left invariant 
metric $\gbar$ on $\Mbar,$ we know that 
the Riemannian volume forms of $(\nilmfldbar1,\gbar)$ and 
$(\nilmfldbar2,\gbar)$ arise from the same left invariant 
volume form on $\Mbar.$  We will denote this volume form 
and its projections onto $(\nilmfldbar1,\gbar)$ and 
$(\nilmfldbar2,\gbar)$ by $\Omega.$
\medskip

Let $\alpha_\G1$ and $\alpha_\G2$ be as in Theorem 3.7. The volume 
forms $\alpha_\G{i}$ are then scalar multiples of $\Omega.$ For $i=1,2,$ 
let $\alpha_\G{i} = p_i\Omega.$
\medskip

So 
$$\int_{\nilmfldbar1}\alpha_\G1 = \ 1 \ = \int_{\nilmfldbar2}\alpha_\G2$$
$$p_1\int_{\nilmfldbar1}\Omega = p_2\int_{\nilmfldbar2}\Omega$$
$$p_1\Vol(\nilmfldbar1,\gbar)= p_2\Vol(\nilmfldbar2,\gbar)$$
\medskip

By hypothesis,  $\spec(\nilmfldbar1,\gbar) = \spec(\nilmfldbar2,\gbar),$ 
and the spectrum of the La\-place-Bel\-tra\-mi operator is known to 
determine the volume of a closed manifold.  Thus 
$\Vol(\nilmfldbar1,\gbar)=\Vol(\nilmfldbar2,\gbar),$ which implies 
$p_1 = p_2,$ and so $\alpha_\G1 = \alpha_\G2.$
\medskip

As the definition of $P_{\alpha_\G{i}}$ depends only on the volume 
form $\alpha_\G{i},$ we must have $$P_{\alpha_\G1}(\tau) = P_{\alpha_\G2}(\tau)$$ 
for all $\tau$ in $\Tau''_1 = \Tau''_2.$ 
Hence $m_1(\tau) = m_2(\tau)$ for all $\tau$ in $\Tau''_1=\Tau''_2.$
\medskip

Thus the representations of $G$
on $\H''_1$ and $\H''_2$ are unitarily equivalent,
and by Proposition 2.2 $$\spec''(\nilmfld1,g) = \spec''(\nilmfld2,g),$$
as desired.
\medskip

The proofs of Lemma 3.4 and Theorem 3.2 are now complete. $\endpf$
\enddemo
\enddemo
\bigskip

\bigskip
\proclaim{Corollary 3.8}
Let $G$ be a simply connected, strictly nonsingular nilpotent Lie group. Two 
cocompact, discrete subgroups {$\G1$ and $\G2$ are representation equivalent 
subgroups of $G$}
if and only if {$\GGGbar1$ and $\GGGbar2$ are representation equivalent 
subgroups of $\Mbar$} and {$\G1 \cap Z(G) = \G2 \cap Z(G).$}  
\endproclaim
\bigskip

\demo{Proof of Corollary 3.8}
\bigskip

Let $\G1$ and $\G2$ be cocompact, discrete subgroups of $G.$
As in the proof of Theorem 3.2, we decompose the representation spaces of 
$\rep1$ and $\rep2$ as 
$$\fnsp{i}= \H'_i \oplus \H''_i,$$ for $i=1,2.$
\medskip

If $\G1$ and $\G2$ are representation equivalent, then the square integrable 
representations occurring in the quasi-regular representations must correspond. 
We showed in the proof of Lemma 3.4  that the square integrable representations 
occurring in $\rep1$ and $\rep2$ are precisely the irreducible representations 
appearing in $\H''_i.$
Thus $\G1$ and $\G2$ are representation equivalent subgroups of $G$ if and only 
if the representations of $G$ on $\H'_1$ and $\H'_2$ are unitarily equivalent 
and the representations of $G$ on $\H''_1$ and $\H''_2$ are unitarily equivalent. 
\medskip

For $i=1,2$ the irreducible components of the representation of $G$ on $\H'_i$ 
correspond to the linear functionals in $\Tau_i$ that are zero on the center of $\mm.$
As these functionals may be viewed as functionals on $\mbar,$ we may likewise 
view the irreducible components as representations of $\Mbar.$  Hence the 
representations of $G$ on $\H'_1$ and $\H'_2$ may be viewed as the quasi-regular 
representations of $\Mbar$ on $\fnspbar1$ and $\fnspbar2,$ respectively.
Thus $\G1$ and $\G2$ are representation equivalent subgroups of $G$ if and only if
$\GGGbar1$ and $\GGGbar2$ are representation equivalent subgroups of $\Mbar$ and the 
representations of $G$ on $\H''_1$ and $\H''_2$ are unitarily equivalent. 
\medskip

Theorem 3.7 tells us that the irreducible representations occurring in 
$\H''_i$ correspond
to the elements of the dual lattice of $\Gamma_i \cap Z(G).$  So if the 
representations of $G$ on $\H''_1$ and $\H''_2$  are unitarily equivalent, 
then the duals of $\Gamma_1\cap Z(G)$ and $\Gamma_2\cap Z(G)$ 
must coincide. Hence  $\Gamma_1\cap Z(G)=\Gamma_2\cap Z(G).$  The proof of the forward
direction is now complete.
\medskip

The reverse direction follows from Proposition 2.2 applied to the quotient 
nilmanifolds and Lemma 3.4.
\endpf
\enddemo 
\bigskip

\bigskip
\bigskip

\subheading{\S4 New Examples of Isospectral Nilmanifolds}
\bigskip

Using Theorem 3.2, we construct and compare five new pairs of isospectral 
nilmanifolds.  A summary of the properties of these examples is listed in 
Section 1, Table I. 
In this Section we compare the quasi-regular representations and the 
fundamental groups of these examples.  The $p$-form spectra of these 
examples are also compared, but the calculations are left to the Appendix.  
The methods used in the Appendix are new, as previously used techniques 
could not be applied to compare the $p$-form spectrum of these higher-step 
examples.   The length spectra and marked length spectra of these examples will be studied in \cite{Gt4}.  
\medskip

With the exception of the column comparing the representation 
equivalence of the fundamental groups,  all of the properties listed in Table I 
are geometric invariants.  Hence a ``No'' in any one of these columns 
demonstrates that an Example is nontrivial.
\medskip

Before proceeding, we need the following.
\bigskip

\proclaim{Definition 4.1} 
Let $\Phi$ be a Lie group automorphism of $G.$
Let $\GG$ be a cocompact, discrete subgroup of $G.$
We say $\Phi$ is an {\it almost inner automorphism} if 
for all elements $x$  of $G$
there exists $a_x$ in $G$ such that   
$\Phi(x)=a_x x a_x^{-1}.$
\endproclaim
\bigskip

\proclaim{Theorem 4.2 (Gordon--Wilson \cite{GW1}) }
Let $G$ be a nilpotent Lie group and let $\GG$ be a cocompact, discrete 
subgroup of $G.$  If $\Phi$ is an almost inner automorphism of $G,$ 
then $\GG$ and $\Phi(\GG)$ are representation equivalent.    \endproclaim
\bigskip

\bigskip
\bigskip

\subheading{Example I}
\bigskip

Consider the simply connected, strictly nonsingular, three-step 
nilpotent Lie group $G$ with Lie algebra
$$\mm = span_{\R}\{X_1, X_2, Y_1, Y_2, Z_1, Z_2, W \}$$
and Lie brackets
$$[X_1,Y_1] = [X_2, Y_2] = Z_1$$
$$[X_1,Y_2] = Z_2$$
$$[X_1,Z_1] = [X_2, Z_2] = [Y_1, Y_2] = W$$
and all other basis brackets zero.
\bigskip

Let $\G1$ be the cocompact, discrete subgroup of $G$ generated by 
$$\{exp(2X_1),exp(2X_2),exp(Y_1),exp(Y_2),exp(Z_1),exp(Z_2),exp(W)\}.$$
Let $\G2$ be the cocompact, discrete subgroup of $G$ generated by 
$$\{exp(2X_1),exp(2X_2),exp(Y_1),exp(Y_2+\frac1{2}Z_2),exp(Z_1),exp(Z_2),exp(W)\}.$$
Note that $\G1 \cap Z(G) = \G2 \cap Z(G) = \{exp(jW):j\in\Z \}.$
\bigskip

Now $\GGGbar2 = \Phi(\GGGbar1)$ where $\Phi$ is the almost inner 
automorphism of $\Mbar$ given on the Lie algebra level by
$$\align \BX_1&\rightarrow \BX_1,\\
\BX_2 &\rightarrow \BX_2,\\
\BY_1 &\rightarrow \BY_1,\\
\BY_2 &\rightarrow \BY_2 + \frac1{2} \BZ_2,\\ 
\BZ_1 &\rightarrow \BZ_1,\\
\BZ_2 &\rightarrow \BZ_2.\endalign$$ 
\bigskip 

By Theorem 4.2, \ $\GGGbar1$ and $\GGGbar2$ are representation equivalent 
subgroups of $\Mbar.$    By Corollary 3.8, \ $\G1$ and $\G2$ are 
representation equivalent subgroups of $G.$ 
\medskip
 
By Proposition 2.2, for any choice of left invariant metric $g$ of $G,$ 
!x
we have 
\break
$p\hy spec(\nilmfld1, g) = p\hy spec(\nilmfld2, g),$ for $p=0,1,\cdots,7.$
\bigskip

\proclaim{Proposition 4.3 (\cite{Gt1})}
The subgroups $\Gamma_1$ and $\Gamma_2$ are not isomorphic as groups.
\endproclaim
\bigskip

\Remark
The author previously established the representation equivalence 
of $\G1$ and $\G2$  in \cite{Gt1} by using a direct calculation.
This example was presented in \cite{Gt1} as the first example of a 
pair of nonisomorphic, representation equivalent subgroups of a solvable 
Lie group. Note that a nilpotent Lie group is necessarily solvable.  
Also, this example was presented in \cite{Gt2} as the first example of a 
pair of representation equivalent subgroups of a solvable Lie group 
producing nilmanifolds with
unequal length spectra.  Contrast this example with what must happen 
in the two-step case.
\bigskip

\proclaim{Definition 4.4}
Let $G$ be a two-step nilpotent Lie group and let $\GG$ be a cocompact, 
discrete subgroup of $G.$ We call the automorphism $\Phi$ of $G$ a 
$\GG${\it -equivalence} if for all $\gamma$ in $\GG$ there exists 
$a_\gamma$ in $G$  and $\gamma'_\gamma$ in $\GG \cap \Gup1$ such that 
$\Phi(\gamma)=a_\gamma \gamma a_\gamma^{-1} \gamma'_\gamma.$
\endproclaim
\bigskip

\proclaim{Theorem 4.5 (\cite{Gt1}, \cite{Gt2})}
Let $G$ be a two-step nilpotent Lie group. Let $\G1$ and $\G2$ be 
cocompact, discrete subgroups of $G$. The subgroups $\G1$ and $\G2$ 
are representation equivalent if and only if there exists $\Phi,$ 
a $\G1$-equivalence of $G,$ such that 
$\Phi(\G1)=\G2.$ Thus if $\G1$ and $\G2$ are representation equivalent, 
they are 
necessarily isomorphic.
In addition, if $\G1$ and $\G2$ are representation equivalent, then $(\nilmfld1,g)$ and
$(\nilmfld2,g)$ have the same length spectrum for any choice of left invariant metric
$g$ of $G.$
\endproclaim
\bigskip
 
\bigskip


\bigskip
\subheading{Example II}
\bigskip

Consider the simply connected, strictly nonsingular, three-step 
nilpotent Lie group $G$ with Lie algebra
$$\mm = span_{\R}\{X_1,Y_1, Y_2, Z,  W \}$$
and Lie brackets
$$[X_1,Y_1] = Z$$
$$[X_1,Z] = [Y_1, Y_2] = W$$
and all other basis brackets zero.
\bigskip

Let $\G1$ be the cocompact, discrete subgroup of $G$ generated by $$\{exp(2X_1),exp(Y_1),exp(Y_2),exp(Z),exp(W)\}.$$
Let $\G2$ be the cocompact, discrete subgroup of $G$ generated by
$$\{exp(2X_1),exp(Y_1+\frac1{2}Z),exp(Y_2),exp(Z),exp(W)\}.$$
Note that $\G1 \cap Z(G) = \G2 \cap Z(G) = \{exp(jW):j\in\Z \}.$
\bigskip

Now $\GGGbar2 = \Phi(\GGGbar1)$ where $\Phi$ is the inner automorphism 
of $\Mbar$ given on the Lie algebra level by 
$$\align \BX_1 &\rightarrow \BX_1, \\
\BY_1 &\rightarrow \BY_1 + \frac1{2} \BZ,\\
\BY_2 &\rightarrow \BY_2, \\
\BZ &\rightarrow \BZ.\endalign$$  Note that an inner automorphism is 
necessarily almost inner.
\medskip
  
By Theorem 4.2, \ $\GGGbar1$ and $\GGGbar2$ are representation equivalent 
subgroups of $\Mbar.$   By Corollary 3.8,  \ $\G1$ and $\G2$ are
representation equivalent subgroups of $G.$
\bigskip

By Proposition 2.2,  $p\hyspec(\nilmfld1, g) = p\hyspec(\nilmfld2, g)$  
for $p=0, 1, \cdots, 5,$
for any choice of left invariant metric $g$ of $G.$ 
\bigskip

Here $\G1$ and $\G2$ are isomorphic.  
Indeed, a simple calculation shows that the isomorphism $\Psi$ given on 
the Lie algebra level by
$$\align
X_1 &\rightarrow X_1 + \frac1{2}Y_2,\\
Y_1 &\rightarrow Y_1 + \frac1{2}Z, \\
Y_2 &\rightarrow Y_2,\\
Z &\rightarrow Z,\\
W &\rightarrow W,
\endalign$$
is an isomorphism of $G$ such that $\Psi(\G1) = \G2.$
\bigskip

\proclaim{Proposition 4.6}
No isomorphism between $\G1$ and $\G2$ will project to a 
$\GGGbar1$-equivalence of $\Mbar.$
\endproclaim
\bigskip

\Remark Example I illustrates that in the higher-step case, the 
representation equivalence of $\G1$ and $\G2$ and the isomorphism 
class of $\G1$ and $\G2$ need not be related. Example II shows that, 
in contrast to Theorem 4.5, even in the case where $\G1$ and $\G2$ 
are isomorphic, knowing the isomorphisms between $\G1$ and $\G2$ 
is not enough to use Corollary 3.8 to establish whether or not 
$\G1$ and $\G2$ are representation equivalent. 
\bigskip

\demo{Proof of Proposition 4.6}
\bigskip

Let $\Psi$ be an isomorphism from $\G1$ to $\G2.$  Extend it to 
the Lie group isomorphism $\Psi:G \rightarrow G$ such that $\Psi(\G1)=\G2.$
\medskip

On the Lie algebra level, any such isomorphism must preserve the 
following ideals of $\mm:$

(1)  $\m2 = span_\R\{W\},$

(2)  $\m1 = span_\R\{Z,W\},$

(3)  $\aa = span_\R\{Y_2,Z,W\},$ and 

(4)  $\cc = span_\R\{Y_1,Y_2,Z,W\},$ 
the centralizer of $\m1$ in $\mm.$

\flushpar To see (3), note that $ad(U)(\mm) \subset \m2$ if and 
only if $U \in \aa.$ \medskip

Note that the generators of $\G1$ and $\G2$ presented above are 
canonical in the sense that every element of $\G1$ may be expressed 
uniquely as $$exp(2n_1X_1)exp(m_1Y_1)exp(m_2Y_2)exp(kZ)exp(jW),$$ 
for integers $n_1,m_1,m_2,k,j.$  Similarly for $\G2.$
\medskip

As $\Psi(\G1) = \G2,$ 
generators of $\G1$ must go to generators of $\G2,$ 
and these generators must be expressible in terms of the canonical 
generators of $\G2,$ given above.
\medskip

Combining this fact with properties (1) through (4), we obtain:
\medskip

$\Psi_*(W) = \pm W,$ \
by (1).

$\Psi_*(Z) = \pm Z + h_0W$ \ 
using (2).

$\Psi_*(Y_2) = \pm Y_2$ mod $\m1,$ \
using (3).

$\Psi_*(Y_1) = \pm (Y_1 + \frac1{2}Z) + h_1 Y_2 + h_2Z$ mod 
$\m2$
using (4).

$\Psi_*(X_1) = \pm X_1 +  \frac1{2}h_3Y_1 + \frac1{2}h_4Y_2$ mod 
$\m1,$
\flushpar where $h_0, \ h_1,\ h_2, \ h_3$ and $h_4$ are integers.
\bigskip

Finally, we use the fact that  $\Psi_*$ is a Lie algebra isomorphism.
By examining the $W$ coefficient of $\Psi_*([X_1,Y_1]) = [\Psi_*(X_1),\Psi_*(Y_1)]$ 
we have the equation $$h_0 = \pm\frac1{2} + \frac1{2}h_1h_3 \pm h_2 \pm\frac1{2}h_4.$$
Thus either $h_3 \neq 0$ or $h_4 \neq 0.$ 
\medskip

As $\BY_1$ and $\BY_2$ are {not} in $[\BX_1,\mbar]$ and {not} in $\mbar^{(1)},$ 
we see that the projection of $\Psi$ cannot possibly be a $\GGGbar1$-equivalence. \endpf 
\enddemo
\bigskip

\bigskip


\bigskip

\subheading{Example III}
\bigskip

Consider again the seven-dimensional Lie group $G$ presented in Example I.
\bigskip 

We again let $\G1$ be the cocompact, discrete subgroup of $G$ generated by
$$\{exp(2X_1),exp(2X_2),exp(Y_1),exp(Y_2),exp(Z_1),exp(Z_2),exp(W)\},$$
and let $\G2$ be the cocompact, discrete subgroup of $G$ generated by
$$\{exp(X_1),exp(X_2), exp(2Y_1),exp(2Y_2), exp(Z_1),exp(Z_2),exp(W)\}.$$
Note that $\G1 \cap Z(G) = \G2 \cap Z(G) = \{exp(jW):j\in\Z\}.$
\bigskip
  
Let $g$ be the left invariant metric on $G$ defined by letting 
$$\{X_1, X_2, Y_1, Y_2, Z_1, Z_2, W\}$$ be an orthonormal basis of $\mm.$  
\bigskip

Now $\GGGbar2 = \Phi(\GGGbar1)$ where $\Phi$ is the automorphism of $\Mbar$ 
given on the Lie algebra level by 
$$\align \BX_1 &\rightarrow \BY_2,  \\
 \BX_2 &\rightarrow \BY_1, \\
 \BY_1 &\rightarrow \BX_2, \\
 \BY_2 &\rightarrow \BX_1, \\
 \BZ_1 &\rightarrow -\BZ_1, \\
 \BZ_2 &\rightarrow -\BZ_2.\endalign$$
\medskip

The automorphism $\Phi$ is also an isometry of $(\Mbar,\gbar),$
and an isometry must preserve the spectrum. Thus 
$\spec(\nilmfldbar1, \gbar) = \spec(\nilmfldbar2, \gbar).$
By Theorem 3.2,  $\spec(\nilmfld1, g) = \spec(\nilmfld2, g).$
\bigskip

\proclaim{Proposition 4.7} The manifolds $(\nilmfld1, g)$ and $(\nilmfld2, g)$ 
are not isospectral on one-forms. 
\endproclaim 
\bigskip

The proof of Proposition 4.7 is left to the Appendix.  Note that 
Proposition 4.7 and Proposition 2.2 together imply that $\G1$ 
and $\G2$ are not representation equivalent subgroups of $G.$
\bigskip

\proclaim{Proposition 4.8}
The subgroups $\G1$ and $\G2$ are not isomorphic as groups.
\endproclaim
\bigskip

\demo{Proof of Proposition 4.8}
\medskip

If there existed a group isomorphism between $\G1$ and $\G2,$
it would extend to a Lie group automorphism $\Psi$ of  $G$
such that $\Psi(\G1)=\G2.$
\medskip

Each of the following ideals of $\mm$ must be preserved by 
the Lie algebra automorphism $\Psi_*:$

(1) $\m2 = \xi = span_{\R}\{ W \},$

(2) $\m1 = span_{\R}\{ Z_1, Z_2, W \},$

(3) $\cc = span_{\R}\{ Y_1, Y_2, Z_1, Z_2, W \},$
the centralizer of $\m1$ in $\mm,$ and

(4) $\aa = span_{\R}\{ X_2, Y_1, Z_1, Z_2, W \}.$

\flushpar To see (4), note that the image of $ad(U)$ has dimension less than three
if and only if $U \in \aa.$ 
\medskip

Now $\Psi(\G1) = \G2.$ Consequently, 
generators of $\G1$ must go to generators of $\G2,$ 
and these generators must be expressible in terms of the canonical 
generators of $\G2,$ given above.
Combining this fact with properties (1) through (4), we obtain:
\medskip
 
$\Psi_*(Y_1) = \pm 2Y_1$ mod $\m1$ by (3) and (4).

$\Psi_*(Y_2) = \pm 2Y_2$ mod $span_\R\{Y_1,Z_1,Z_2,W\}$  by (3).

$\Psi_*(W) = \pm W$ by (1).

\flushpar But then $\Psi_*([Y_1,Y_2]) \neq [\Psi_*(Y_1),\Psi_*(Y_2)].$
\endpf
\enddemo
\bigskip

\bigskip


\bigskip

\subheading{Example IV}
\bigskip

Consider again the five-dimensional Lie group $G$ presented in Example II.
\bigskip 

We again let $\G1$ be the cocompact, discrete subgroup of $G$ generated by
$$\{exp(2X_1),exp(Y_1),exp(Y_2),exp(Z),exp(W)\},$$
and let $\G2$ be the cocompact, discrete subgroup of $G$ generated by 
$$\{exp(X_1), exp(2Y_1), exp(Y_2),exp(Z),exp(W)\}.$$
Note that $\G1 \cap Z(G) = \G2 \cap Z(G) = \{exp(jW):j\in\Z\}.$
\bigskip
  
Let $g$ be the left invariant metric on $G$ defined by letting 
$$\{X_1, Y_1, Y_2, Z, W\}$$ be an orthonormal basis of $\mm.$  
\bigskip

Now $\GGGbar2 = \Phi(\GGGbar1)$ where $\Phi$ is the automorphism of 
$\Mbar$ given on the Lie alegebra level by 
$$\align \BX_1 &\rightarrow \BY_1,\\
\BY_1 &\rightarrow \BX_1,\\
\BY_2 &\rightarrow \BY_2, \\
\BZ &\rightarrow -\BZ.\endalign$$
\medskip

The automorphism $\Phi$ is clearly an isometry of $(\Mbar,\gbar),$
and an isometry preserves the spectrum. Thus 
$\spec(\nilmfldbar1, \gbar) = \spec(\nilmfldbar2, \gbar).$
By Theorem 3.2, $\spec(\nilmfld1, g) = \spec(\nilmfld2, g).$
\bigskip

\proclaim{Proposition 4.9}
The manifolds $(\nilmfld1, g)$ and $(\nilmfld2, g)$ are not isospectral on one-forms. 
\endproclaim
\bigskip

The proof of Proposition 4.9 is left to the Appendix.
\bigskip

\proclaim{Proposition 4.10}
The subgroups $\G1$ and $\G2$ are not isomorphic as groups.
\endproclaim
\bigskip

\Remark The combination of properties exhibited by Examples III and IV 
are similar to properties exhibited by pairs of isospectral Heisenberg 
manifolds constructed by Gordon and Wilson \cite{GW2, G2}.
\bigskip

\demo{Proof of Proposition 4.10}
\medskip

If there existed a group isomorphism between $\G1$ and $\G2,$
it would extend to a Lie group automorphism $\Psi$ of  $G$
such that $\Psi(\G1)=\G2.$
\medskip

As before, generators of $\G1$ must go to generators of $\G2,$ 
and these generators must be expressible in terms of the canonical 
generators of $\G2,$ given above.
Combining this with properties (1) through (4) from the proof of 
Proposition 4.6, we have

$\Psi_*(Y_1) = \pm 2Y_1$ mod $span_\R\{Y_2,Z,W\}.$

$\Psi_*(Y_2) = \pm Y_2$ mod $\m1.$

$\Psi_*(W) = \pm W.$ 

\flushpar But then $\Psi_*([Y_1,Y_2]) \neq [\Psi_*(Y_1),\Psi_*(Y_2)].$
\endpf
\enddemo
\bigskip

\bigskip


\bigskip

\subheading{Example V}
\bigskip

Consider again the seven-dimensional Lie group $G$ presented in Example I.
\bigskip 

We fix a left invariant metric on $G$ by letting 
$\{E_1,E_2,E_3,E_4,E_5,E_6,E_7\}$ be an orthonormal basis of $\gg$ where
$$
\align
E_1&= X_1-\frac{1}{2}X_2-\frac{1}{4}Y_2,\\
E_2&= X_2-\frac{1}{4}Y_1,\\
E_3&= Y_1,\\
E_4&=Y_1+Y_2,\\
E_5&=Z_1,\\
E_6&=\frac{1}{2}Z_1+Z_2,\\
E_7&=W.
\endalign
$$

Let $\Phi$ be the automorphism of $G$ defined on the Lie algebra level by 
$$
\align
X_1&\rightarrow -X_1+X_2+\frac{1}{4}Y_1+\frac{1}{2}Y_2,\\
X_2&\rightarrow X_2-\frac{1}{2}Y_1+\frac{1}{4}Z_1,\\
Y_1&\rightarrow -Y_1,\\
Y_2&\rightarrow 2Y_1+Y_2+Z_2,\\
Z_1&\rightarrow Z_1+\frac1{2}W,\\
Z_2&\rightarrow -Z_1-Z_2+\frac1{4}W,\\
W&\rightarrow -W.
\endalign
$$
A straightforward calculation shows that $\Phi_*([U,V])=[\Phi_*(U), \Phi_*(V)]$ 
for all $U,V$ in $\gg.$  Thus $\Phi$ is indeed a Lie group automorphism.
\bigskip

Let $\G1$ be the cocompact, discrete subgroup of $G$ generated by $$\{exp(2X_1),exp(2X_2),exp(Y_1),exp(Y_2),exp(Z_1),exp(Z_2),exp(W)\},$$
and let $\G2=\Phi(\G1).$  Note that $\G1 \cap Z(G) = \G2 \cap Z(G) = \{exp(jW):j\in\Z\}.$  
\bigskip

Let $\bar\Phi$ be the projection of $\Phi$ onto $\Gbar.$
Then $\bar\Phi$ factors as $\bar\Phi = \Psi_1 \circ \Psi_2$ where $\Psi_1$ is the 
automorphism of $\Mbar$ given on the Lie algebra level by 
$$\align \BX_1 &\rightarrow -\BX_1+\BX_2+\frac{1}{4}\BY_1+\frac{1}{2}\BY_2, \\
 \BX_2 &\rightarrow \BX_2-\frac{1}{2}\BY_1,\\
 \BY_1 &\rightarrow -\BY_1,\\
 \BY_2 &\rightarrow 2\BY_1+\BY_2,\\
 \BZ_1 &\rightarrow \BZ_1,\\
 \BZ_2 &\rightarrow -\BZ_1-\BZ_2,
\endalign$$
and $\Psi_2$ is the automorphism of $\Mbar$ given on the Lie algebra level by 
$$\align \BX_1 &\rightarrow \BX_1,\\
 \BX_2 &\rightarrow \BX_2+\frac{1}{4}\BZ_1,\\
 \BY_1 &\rightarrow \BY_1,\\
 \BY_2 &\rightarrow \BY_2-\BZ_1-\BZ_2,\\
 \BZ_1 &\rightarrow \BZ_1,\\
 \BZ_2 &\rightarrow \BZ_2.
\endalign$$
\medskip

By rewriting $\Psi_1$ in terms of the orthonormal basis 
$\{\BE_1, \BE_2, \BE_3, \BE_4, \BE_5, \BE_6 \}$ of $\gggbar,$
one easily sees that $\Psi_1(\BE_i)=\pm\BE_i$ for $i=1, \dots, 6.$  Thus the 
automorphism $\Psi_1$ is also an isometry of $\Gbar$ and must preserve the spectrum.
A simple calculation shows that $\Psi_2$ is an almost inner automorphism of $\Gbar,$ 
which by Theorem 4.2 also preserves the spectrum. 
Thus $\spec(\nilmfldbar1, \gbar) = \spec(\nilmfldbar2, \gbar).$
By Theorem 3.2, $\spec(\nilmfld1, g) = \spec(\nilmfld2, g).$
\bigskip
 
\proclaim{Proposition 4.11}
The manifolds $(\nilmfld1, g)$ and $(\nilmfld2, g)$ are not isospectral on one-forms. 
\endproclaim
\bigskip

The proof of Proposition 4.11 is left to the Appendix.
\bigskip

\Remark We will show in \cite{Gt4} that the automorphism $\Phi$ marks the 
length spectrum of these examples.  This is the first example of a pair of 
manifolds with the same marked length spectrum but not the same spectrum on one-forms.
\bigskip

\bigskip


\bigskip

\subheading{Appendix:  Comparing the $p$-form spectrum of nilmanifolds}
\bigskip

In this appendix, we show that the pairs of isospectral manifolds in 
Examples III, IV, and V are not isospectral on one-forms.  
\bigskip

Recall that on smooth $p$-forms, the Laplace-Beltrami operator is defined as
$$\Delta = d \delta + \delta d.$$ 
Here $\delta$ is the metric adjoint of $d.$  
Equivalently, $\delta = (-1)^{n(p+1)+1}*d*$ where $*$ is the Hodge-$*$ operator.
Let $E^p(M)$ denote the exterior algebra of smooth differential $p$-forms on $M.$	
Then for $f \in C^\infty(M)$ and $\tau \in E^p(M),$ Gordon 
and Wilson \cite {GW1} showed that 
$$\Delta(f \tau) = (\Delta f)\tau + f(\Delta \tau) - 2\nabla_{\text{\it grad } f}\tau.\tag{A.1}$$
\bigskip

For $G$ a simply connected Lie group with cocompact, 
discrete subgroup $\GG,$ 
view $E^p(\GG \backslash G)$ as 
$$E^p(\GG \backslash G) = C^\infty(\GG \backslash G) \otimes \Lambda^p(\gg^*).$$ 
Here elements of $\Lambda^p(\gg^*)$ are viewed both as left 
invariant $p$-forms of $G$ and also as elements of $E^p(\GG \backslash G).$
\bigskip

\proclaim{Proposition 4.7}
The nilmanifolds $(\nilmfld1, g)$ and $(\nilmfld2,g)$ as presented in 
Example III are not isospectral on one-forms.  
\endproclaim
\bigskip

\demo{Outline of Proof}
\medskip

Step 1:  We decompose $\one  \spec(\nilmfld{i}, g)$ into four
components:  
$$\one \spec^{I}(\nilmfld{i}, g) \cup \one \spec^{II}(\nilmfld{i}, g) 
\cup \one \spec^{III}(\nilmfld{i}, g) \cup \one \spec^{IV}(\nilmfld{i},g).$$
The multiplicity of an eigenvalue in $\one \spec(\nilmfld{i},g)$ is 
the sum of its multiplicities in each of the four components.
\bigskip

Step 2:  Using representation theory, we show 
$$\one \spec^{IV}(\nilmfld1, g)=\one \spec^{IV}(\nilmfld2, g)$$ 
and $$\one \spec^{III}(\nilmfld1,g) = \one \spec^{III}(\nilmfld2,g).$$
\bigskip

Step 3:  We show that the multiplicity of every eigenvalue 
in $\one \spec^{II}(\nilmfld{i}, g)$ is congruent to $0$ modulo $4.$
\bigskip

Step 4:  Finally, we show that the eigenvalue $\pi^2 + 1$ 
does not occur in 
$\one \spec^I(\nilmfld1, g)$ but occurs with multiplicity 2 
in $\one \spec^I(\nilmfld2, g).$ 
\bigskip

The result now follows. $\endpf$ 
\enddemo
\bigskip

\demo{Proof of Proposition 4.7}

\demo{Step 1}
\bigskip

Using the notation of Section 2, for 
$i=1,2,$ let $\Tau_{i}$ be a subset of $\mm^*$ such that 
$$\fnsp{i} \cong {\underset{\tau\in\Tau_{i}}\to\bigoplus} {m_i(\tau)}\H_\tau.$$
\medskip 

Let $\Tau_i = \Tau^{I}_i \cup \Tau^{II}_i \cup \Tau^{III}_i \cup \Tau^{IV}_i$
where
$$\align \Tau^{I}_i &= \{\tau \in \Tau_i :  \tau(Z_1) = \tau(Z_2) = \tau(W) = 0\}, \\
\Tau^{II}_i &= \{\tau \in \Tau_i :  \tau(Z_2) \neq 0, \tau(Z_1) = \tau(W) = 0 \},\\
\Tau^{III}_i &= \{\tau \in \Tau_i :  \tau(Z_1) \neq 0, \tau(W) = 0\},\\
\Tau^{IV}_i &= \{\tau \in \Tau_i : \tau(W) \neq 0 \}.
\endalign$$
\medskip

Let
$$\align
\H^{I}_i ={\underset{\tau\in\Tau^{I}_{i}}\to\bigoplus}{m_i(\tau)}\H_\tau ,
&\qquad
\H^{II}_i ={\underset{\tau\in\Tau^{II}_{i}}\to\bigoplus}{m_i(\tau)}\H_\tau,\\
\H^{III}_i ={\underset{\tau\in\Tau^{III}_{i}}\to\bigoplus}{m_i(\tau)}\H_\tau,
&\qquad 
\H^{IV}_i ={\underset{\tau\in\Tau^{IV}_{i}}\to\bigoplus}{m_i(\tau)}\H_\tau.
\endalign$$
\medskip

As representation spaces, 
$$\fnsp{i} = \H^{I}_i \oplus \H^{II}_i \oplus \H^{III}_i \oplus \H^{IV}_i.$$
\medskip

Decompose $\one \spec(\nilmfld{i},g)$ as
$$\one \spec^{I}(\nilmfld{i}, g) \cup \one \spec^{II}(\nilmfld{i}, g) 
\cup \one \spec^{III}(\nilmfld{i}, g) \cup \one \spec^{IV}(\nilmfld{i},g),$$ 
where 
$\one \spec^I(\nilmfld{i}, g)$ 
is defined as the spectrum of the Laplacian 
acting on 
$\H^I_i \otimes \Lambda^1(\mm^*).$ Define $\one \spec^{II}(\nilmfld{i},g), \
\one \spec^{III}(\nilmfld{i},g),$ and 
$\one \spec^{IV}(\nilmfld{i},g)$ similarly.  The multiplicity of an eigenvalue in 
$\one \spec(\nilmfld{i}, g)$ 
is equal to the sum of its multiplicities in each of the four components.
\enddemo
\bigskip 

\bigskip


\demo{Step 2}
\medskip

By Lemma 3.4,  the representations of 
$G$ on $\H^{IV}_1$ and $\H^{IV}_2$ are unitarily equivalent and 
$$\one\spec^{IV}(\nilmfld1, g) =  \one\spec^{IV}(\nilmfld2, g).$$
\medskip

We now show that the representations of $G$ on $\H^{III}_1$
and $\H^{III}_2$ are unitarily equivalent.
\medskip

The irreducible representations of $G$ corresponding to 
elements of $\mm^*$ that are zero on the center may be viewed
as irreducible representations of $\Mbar = G / Z(G).$
It is easy to see that such representations of $G$
are unitarily equivalent if and only if
the corresponding representations of $\Mbar$
are unitarily equivalent.
\medskip

For all $\tau \in \Tau^{III}_i, \tau(\xi)=0.$  
We may thus use the following Proposition to calculate 
the representations of $G$ on $\H^{III}_1$ and $\H^{III}_2.$
\bigskip

\proclaim{Proposition A.2 (see Pesce \cite{P2})}
Let 
$N$ 
be a simply connected, two-step nilpotent Lie group with 
$\GG$ a cocompact, discrete subgroup of 
$N.$  Let 
$\tau \in \nn^*.$   Define 
$\nn_\tau = \{ Y \in \nn : \tau([Y,\nn]) \equiv 0 \}.$  
Then the irreducible representation 
$\pi_\tau$ appears in the quasi-regular representation of 
$N$ on $L^2(\GG \backslash N)$ 
if and only if 
$\tau(log \GG \cap \nn_\tau) \subset \Z.$  
The multiplicity of 
$\pi_\tau$ is one if 
$\tau(\nn^{(1)}) \equiv 0.$   
Define a nondegenerate, skew-symmetric bilinear form on 
$\nn / \nn_\tau$ by 
$B_\tau(U,V)=\tau([U,V])$ 
for all $U,V \in \nn / \nn_\tau.$  
Then if 
$\tau(\nn^{(1)}) \not \equiv 0,$ 
the multiplicity of 
$\pi_\tau$ is equal to 
$\sqrt{det B_\tau},$  
where the determinant is calculated with respect to any basis of 
${\Cal L}_\tau = log \GG /( log \GG \cap \nn_\tau ).$
\endproclaim
\bigskip

\Remark The occurrence condition above actually follows directly 
from a more general occurrence and multiplicity theorem due 
independently to Richardson \cite {R} and Howe \cite{H}.  
\bigskip

Let 
$\{\alpha_1, \alpha_2, \beta_1, \beta_2, \zeta_1, \zeta_2, \omega\}$ 
be the dual basis to the orthonormal basis 
$\{X_1, X_2, Y_1, Y_2, \allowbreak Z_1, \allowbreak Z_2, W\}$ of 
$\mm.$
\medskip

If 
$\tau \in \Tau^{III}_i,$  then 
$$\tau = A_1\alpha_1 + A_2 \alpha_2 + B_1 \beta_1 + B_2 \beta_2  + C_1 \zeta_1 + C_2 \zeta_2$$
for some $A_1, A_2, B_1, B_2, C_1, C_2 \in \R$  with $C_1 \neq 0.$  
\medskip

Now 
$\mbar_\tau = span_\R\{\BZ_1,\BZ_2\}.$
Hence 
$$\split
log\GGGbar1 \cap \mbar_\tau &= log\GGGbar2 \cap \mbar_\tau \\
& = log(exp(\Z \BZ_1)exp(\Z \BZ_2))\\
& = span_\Z\{\BZ_1,\BZ_2\}. \endsplit$$   
So $\tau(log\GGGbar{i} \cap \mbar_\tau)\subset \Z$ 
if and only if
$C_1 \in \Z$ and $C_2 \in \Z.$  
By Proposition A.2, we see that
$\tau \in \Tau^{III}_i$  
if and only if $C_1 \in \Z$ and $C_2 \in \Z.$  
Moreover, distinct values of $C_1$ and $C_2$ determine distinct 
coadjoint orbits of $\mm^*.$  As these conditions are the same for both $\Tau^{III}_1$ and $\Tau^{III}_2,$ 
we may assume $\Tau^{III}_1 = \Tau^{III}_2.$
\medskip

We now calculate multiplicities.
\medskip

A basis for 
${\Cal L_1}_\tau = log\GGGbar1 / (log\GGGbar1 \cap \mbar_\tau)$
is $\{2\BX_1, 2\BX_2, \BY_1, \BY_2\}.$  
A basis for 
\break
${\Cal L_2}_\tau = log\GGGbar2 / (log\GGGbar2 \cap \mbar_\tau)$
is $\{\BX_1, \BX_2, 2\BY_1, 2\BY_2\}.$  
In both cases, $\sqrt{det B_\tau}= 4{C_1}^2.$
Thus for $\tau \in \Tau^{III}_1 = \Tau^{III}_2,$ the multiplicities $m_1(\tau)$ 
and $m_2(\tau)$ are equal. Hence the representations of $\Mbar$ are unitarily 
equivalent, so the representations of $G$ on $\H^{III}_1$ and $\H^{III}_2$ are 
unitarily equivalent.
\bigskip

By Proposition 2.2,  
$$\one spec^{III}(\nilmfld1, g) = \one spec^{III}(\nilmfld2, g),$$
as desired. 
\enddemo
\bigskip

\bigskip

\demo{Step 3}
\medskip

We now show that for any eigenvalue 
in $\one spec^{II}(\nilmfld{i}, g),$ its multiplicity in 
\break
$\one spec^{II}(\nilmfld{i}, g)$ is always congruent to 0 modulo 4.
\medskip

We first compute the multiplicity of the irreducible representations 
occurring here, using the same technique as in Step 2.
\medskip

For all 
$\tau \in \Tau^{II}_i,$  
$\tau = A_1\alpha_1 + A_2 \alpha_2 + B_1 \beta_1 + B_2 \beta_2  + C_2 \zeta_2$
for some $A_1, \ab A_2, \ab B_1, \ab B_2, \ab C_2\in\R$
with $C_2 \neq 0.$ We may again use Proposition A.2. 
\medskip

Now 
$\mbar_\tau = span_\R\{\BX_2,\BY_1, \BZ_1,\BZ_2\}.$
Hence 
$$\split
log\GGGbar1 \cap \mbar_\tau &= log(exp(2\Z \BX_2)exp(\Z \BY_1)exp(\Z \BZ_1)exp(\Z \BZ_2))\\
& = span_\Z\{2\BX_2, \BY_1, \BZ_1,\BZ_2\}. \endsplit$$   
So $\tau(log\GGGbar1 \cap \mbar_\tau)\subset \Z$ 
if and only if
$A_2 \in \frac1{2}\Z, \ B_1 \in \Z$ and $C_2 \in \Z.$
\medskip

However, 
$$\split
log\GGGbar2 \cap \mbar_\tau &= log(exp(\Z \BX_2)exp(2\Z \BY_1)exp(\Z \BZ_1)exp(\Z \BZ_2))\\
& = span_\Z\{\BX_2, 2\BY_1, \BZ_1,\BZ_2\}. \endsplit$$   
So $\tau(log\GGGbar2 \cap \mbar_\tau)\subset \Z$ 
if and only if
$A_2 \in \Z, \ B_1 \in \frac1{2}\Z$ and $C_2 \in \Z.$
\medskip

We now calculate $m_i(\tau).$
A basis for 
${\Cal L_1}_\tau = log\GGGbar1 / (log\GGGbar1 \cap \mbar_\tau)$
is $\{2\BX_1, \BY_2\}.$  
A basis for 
${\Cal L_2}_\tau 
= log\GGGbar2 / (log\GGGbar2 \cap \mbar_\tau)$
is $\{\BX_1, 2\BY_2\}.$ 
In both cases $\sqrt{det B_\tau}= 2|C_2|.$
\medskip
 
Thus if $\tau$ is in $\Tau^{II}_i,$ the irreducible representation 
$\pi_\tau$ occurs in the representation of $G$ on $\H^{II}_i$ with 
multiplicity $2|C_2|.$ As the integer $C_2 \neq 0,$ the multiplicity 
must be even.  Hence any eigenvalue of $\Delta$ acting on 
$\H^{II}_i \otimes \Lambda^1(\mm^*)$ must occur in 
$\one spec^{II}(\nilmfld{i},g)$ with multiplicity congruent to 0 modulo $2.$  \medskip

We now use the following.
\medskip

\proclaim{Proposition A.3 (Gordon--Wilson, \cite{GW1})}
Let $G$ be a simply connected Lie group with left invariant metric $g,$
and let $\G1$ and $\G2$ be cocompact, discrete subgroups of $G.$
Let $\H_1$ and $\H_2$ be invariant subspaces of $\rep1$ and $\rep2,$ 
respectively.  Denote by $p\hy\spec'(\nilmfld{i},g)$ the spectrum of $\Delta$ 
restricted to acting on $\H_i\otimes\Lambda^p(\gg^*).$ Assume there exists an 
automorphism $\Phi$ of $G$ such that
\flushpar (1)  $\Phi$ is also an isometry of $(G,g),$ \ and
\flushpar (2)  $\rep1$ restricted to $\H_1$ and $\rep2\circ\Phi$
restricted to $\H_2$ are unitarily equivalent,
then
$$p\hy\spec'(\nilmfld1,g)= p\hy\spec'(\nilmfld2,g).$$
Here $\rep2\circ\Phi$ is defined by
$((\rep2\circ\Phi)(x))f=f\circ R_{\Phi(x)}$ for $x$ in $G$ and $f$ in $\fnsp1.$  
\endproclaim
\bigskip 

Let $\Phi$ be the Lie group automorphism of $G$ defined on the Lie algebra level by

$$\align
X_1 & \rightarrow -X_1\\
X_2 & \rightarrow X_2\\
Y_1 & \rightarrow -Y_1\\
Y_2 & \rightarrow Y_2\\
Z_1 & \rightarrow Z_1\\
Z_2 & \rightarrow -Z_2\\
W & \rightarrow -W\endalign$$
\medskip

The automorphism $\Phi$ is also an isometry of $(G,g).$  
Note that 
$\tau \circ \Phi_* = -A_1\alpha_1 + A_2\alpha_2 - B_1\beta_1 + B_2\beta_2 -C_2\zeta_2.$  
Clearly, if $\tau$ satisfies Condition ($*$) or ($**$), then so does $\tau \circ \Phi_*.$    
A straightforward calculation shows that since $C_2 \neq 0,$ the functionals $\tau$ 
and $\tau \circ \Phi_*$ are not in the same coadjoint orbit of $\mm^*,$ 
so $\pi_\tau$ and $\pi_{\tau \circ \Phi_*}$ are not unitarily equivalent.  
Thus if $\pi_\tau$ occurs in $\H^{II}_i$ with multiplicity $2|C_2|$ then 
so does $\pi_{\tau \circ \Phi_*},$ also with multiplicity $2|C_2|.$
\medskip

Note that  by (2.3) $\pi_{\tau \circ \Phi_*} = \pi_\tau \circ \Phi.$ Using 
Proposition A.3,  any eigenvalue of $\Delta$ acting on 
$\H_\tau \otimes \Lambda^1(\mm^*)$ must also occur as an eigenvalue  of 
$\Delta $ acting on 
$\H_{\tau \circ \Phi_*} \otimes \Lambda^1(\mm^*).$ And each of the 
representation spaces $\H_\tau$ and $\H_{\tau\circ\Phi_*}$ occurs
in $\H^{II}_i$ with multiplicity $2|C_2|.$  Consequently,  the 
multiplicity of any eigenvalue in $spec^{II}(\nilmfld{i},g)$ is a 
multiple of $4|C_2|,$  which is clearly congruent to 0 modulo 4, as desired.
\enddemo
\bigskip

\bigskip


\demo{Step 4}
\medskip

We now show that the eigenvalue $\pi^2 + 1$ does not occur 
in $\one spec^I(\nilmfld1, g)$ but occurs with multiplicity 2 
in $\one spec^I(\nilmfld2, g).$
\medskip

For $\tau \in \Tau^I_1$ or $\tau \in \Tau^I_2,\  \tau(\g1) \equiv 0.$  
We again use Proposition A.2 to calculate the irreducible 
representations occurring here.  We write 
$\tau = A_1\alpha_1 + A_2\alpha_2 +  B_1\beta_1 + B_2\beta_2$ for 
some $A_1, A_2, B_1, B_2 \in \R.$
\medskip

Now $\mbar_\tau = \mbar,$ so 
$\tau(log\GGGbar{i} \cap \mbar_\tau) \subset \Z$ if and only if 
$\tau(log \GGGbar{i}) \subset \Z.$  
\medskip

Thus, $\tau \in \Tau^I_1$ if and only if 
$$A_1,A_2 \in \frac1{2}\Z \text{ and } B_1,B_2 \in \Z,\tag{$*$}$$ and
$\tau \in \Tau^I_2$ if and only if 
$$A_1,A_2 \in \Z \text{ and } B_1,B_2 \in \frac1{2}\Z.\tag{$**$}$$ 
\medskip

Let $\H_\tau$ be the associated representation space of 
$\pi_\tau.$  Then  $\H_\tau$ may be viewed as the 
one-dimensional subspace of $\fnsp{i}$ generated by  (see Section 2):
$$\split \F(exp(x_1X_1)&exp(x_2X_2)exp(y_1Y_1)exp(y_2Y_2)exp(z_1Z_1)exp(z_2Z_2)exp(wW)) \\
&= exp\{2\pi i \tau(x_1X_1 + x_2X_2 + y_1Y_1 + y_2Y_2) \}.\endsplit$$
That is, $\H_\tau = \Cplx \F.$
\medskip

We now calculate $\Delta$ acting on $\H_\tau \otimes \Lambda^1(\mm^*).$
Note that if we let $\Cplx \Lambda^1(\mm^*)$ denote $\Lambda^1(\mm^*)$ 
with complex coefficients, then 
$$\H_\tau \otimes \Lambda^1(\mm^*) = \F \otimes \Cplx \Lambda^1(\mm^*).$$
\medskip

Let $\F \otimes \mu \in \F \otimes \Cplx \Lambda^1(\mm^*).$
Then $\mu = a_1\alpha_1 + a_2\alpha_2 + b_1\beta_1 + b_2\beta_2 + z_1\zeta_1 + z_2\zeta_2 + w\omega$ 
for some $a_1, a_2, b_1, b_2, z_1, z_2, w \in \Cplx.$
\medskip

As $\F$ is independent of $z_1, z_2,$ and $w$ we have
$$
\align \Delta \F &= - X_1^2\F - X_2^2\F - Y_1^2\F - Y_2^2\F\\
&= 4\pi^2(A_1^2 + A_2^2 + B_1^2 + B_2^2)\F\\
&= 4\pi^2S^2\F
\endalign$$
where $S^2 = A_1^2 + A_2^2 + B_1^2 + B_2^2.$
\medskip

Let $*$ denote the Hodge-$*$ operator. One easily sees that 
$d*\mu = 0$ for all $\mu \in \mm^*.$  Hence $\delta\mu = \pm*d*\mu = 0$ 
for all $\mu \in \mm^*.$  
Consequently $\Delta = \delta d$ on $\Lambda^1(\mm^*).$  
\medskip

For $\mu\in \Lambda^1(\mm^*),$  $d\mu(U,V)=-\mu([U,V])$ for all 
$U,V \in \mm.$
Using this fact together with the definition of $\delta$ as the 
metric adjoint of $d,$ one
easily computes that  
$\Delta\alpha_1 =\Delta\alpha_2 = \Delta\beta_1 = \Delta\beta_2 =0,$
$\Delta \zeta_1 = 2\zeta_1, \Delta\zeta_2 = \zeta_2,$ and 
$\Delta \omega = 3\omega.$
\medskip

To calculate $\Delta (\F \otimes \mu),$
it remains to calculate  
$$\nabla_{grad\F}\mu = 2\pi i\F(A_1\nabla_{X_1}\mu + A_2\nabla_{X_2}\mu  + B_1\nabla_{Y_1}\mu + B_2\nabla_{Y_2}\mu).$$
\medskip

For Lie algebras with a left invariant metric, the covariant 
derivatives can be calculated via the following equation: 
$$<\nabla_UV,U'> =\frac1{2} <[U',U],V> + \frac1{2}<[U',V],U> + \frac1{2}<[U,V],U'>,$$ 
for $U,V,U'$ left invariant vector fields of $\mm.$
A simple calculation shows that $\nabla_U(V^\flat) = (\nabla_UV)^\flat,$
where $V^\flat$ denotes the dual of $V$ in $\mm^*$ with respect to our choice of 
orthonormal basis; that is, 
$V^\flat(U) =  <V,U>$ for all $U$ in $\mm^*.$
\medskip

We thus obtain the following chart:

{\eightpoint
$$
\vbox{\offinterlineskip
\halign{\strut\vrule#&\quad#\hfil\quad&&\vrule#&\quad\hfil#\hfil\quad\cr
\noalign{\hrule}
&&&&&&&&&&&&&&&&\cr
&$\nabla_U\mu$&&$\alpha_1$&&$\alpha_2$&&$\beta_1$&&$\beta_2$&&$\zeta_1$&&$\zeta_2$&&$\omega$&\cr
&&&&&&&&&&&&&&&&\cr
\noalign{\hrule}
\noalign{\hrule}
&&&&&&&&&&&&&&&&\cr
&$X_1$&&0&&0&&$\frac1{2}\zeta_1$&&$\frac1{2}\zeta_2$&&$-\frac1{2}\beta_1+\frac1{2}\omega$&&$-\frac1{2}\beta_2$&&$-\frac1{2}\zeta_1$&\cr
&&&&&&&&&&&&&&&&\cr
\noalign{\hrule}
&&&&&&&&&&&&&&&&\cr
&$X_2$&&0&&0&&0&&$\frac1{2}\zeta_1$&&$-\frac1{2}\beta_2$&&$\frac1{2}\omega$&&$-\frac1{2}\zeta_2$&\cr
&&&&&&&&&&&&&&&&\cr
\noalign{\hrule}
&&&&&&&&&&&&&&&&\cr
&$Y_1$&&$-\frac1{2}\zeta_1$&&0&&0&&$\frac1{2}\omega$&&$\frac1{2}\alpha_1$&&0&&$-\frac1{2}\beta_2$&\cr
&&&&&&&&&&&&&&&&\cr
\noalign{\hrule}
&&&&&&&&&&&&&&&&\cr
&$Y_2$&&$-\frac1{2}\zeta_2$&&$-\frac1{2}\zeta_1$&&$-\frac1{2}\omega$&&0&&$\frac1{2}\alpha_2$&&$\frac1{2}\alpha_1$&&$\frac1{2}\beta_1$&\cr
&&&&&&&&&&&&&&&&\cr
\noalign{\hrule}
}}
$$
}
\bigskip 

Using (A.1) and the above information, a straightforward calculation shows that 
if we let $E_\tau =$

$$\left(\matrix 
4\pi^2S^2 & 0 & 0 & 0 &  -2\pi i B_1 & -2\pi i B_2 & 0\\
0 & 4 \pi^2S^2 & 0 & 0 & -2 \pi i B_2 & 0 & 0\\
0 & 0 & 4 \pi^2S^2 & 0 & 2 \pi i A_1 & 0   & -2 \pi i B_2\\
0 & 0 & 0 & 4 \pi^2S^2 & 2 \pi i A_2& 2 \pi i A_1 & 2\pi i B_1\\
2\pi i B_1 &2 \pi i B_2 & -2 \pi i A_1 & -2 \pi i A_2 &  4 \pi ^2 S^2 + 2 &0 & 2 \pi i A_1 \\
2\pi i B_2 & 0 & 0 & -2 \pi i A_1 & 0 &  4 \pi ^2 S^2 + 1 & 2 \pi i A_2 \\
0 & 0 & 2 \pi i B_2 & -2 \pi i B_1 & -2 \pi i A_1 & -2 \pi i A_2 & 4 \pi^2S^2 + 3
\endmatrix
\right),$$
then $\Delta(\F \otimes \mu) = \lambda(\F \otimes \mu)$ if and only if 
$\lambda$ is an eigenvalue of the matrix  $E_\tau.$
\medskip

We now calculate necessary conditions on 
$\tau=A_1\alpha_1+A_2\alpha_2+B_1\beta_1+B_2\beta_2$ so that 
$\pi^2 + 1$ is an eigenvalue of $E_\tau.$
As $\tau \in \Tau^I_1$ or $\tau \in \Tau^I_2,$  we know  
$A_1, A_2,B_1, B_2 \in \Q.$
If $det(E_\tau - (\pi^2+1)I_7) = 0,$ then $\pi$ is the root of a 
polynomial with rational coefficients.  However $\pi$ is transcendental.  
Thus the coefficients of the powers of $\pi$ must be zero.  
\medskip

A straightforward calculation shows that $\pi^{14}$ is the highest power of 
$\pi$ occurring in the polynomial, and the coefficient of $\pi^{14}$ 
is equal to $(4S^2-1)^7.$  Thus if $\pi^2+1$ is an eigenvalue of 
$E_\tau,$ then $S^2 = \frac1{4}.$  Recall $S^2 = A_1^2 + A_2^2 + B_1^2 + B_2^2.$
\medskip

For $\tau$ in $\Tau^I_1,$ \ $S^2=\frac{1}{4}$ if and only if (see ($*$))  
$\tau = \pm \frac1{2}\alpha_1$ or $\tau = \pm \frac1{2}\alpha_2.$
And for $\tau$ in $\Tau^I_2,$ \  $S^2=\frac{1}{4}$
if and only if (see ($**$)) $\tau = \pm \frac1{2}\beta_1$ or 
$\tau = \pm \frac1{2}\beta_2.$
\medskip

For $\tau = \pm\frac1{2}\alpha_1, \tau =\pm\frac1{2}\alpha_2,$ or 
$\tau = \pm\frac1{2}\beta_2,$
a simple calculation shows that 
$det(E_\tau - (\pi^2+1)I_7) \neq 0.$  Thus, the eigenvalue $\pi^2 + 1$ 
does not arise from the Laplacian acting on 
$\H_{\pm\frac1{2}\alpha_1} \otimes \Lambda^1(\mm^*),$ 
$\H_{\pm\frac1{2}\alpha_2} \otimes \Lambda^1(\mm^*),$ or 
$\H_{\pm\frac1{2}\beta_2} \otimes \Lambda^1(\mm^*),$
and $$\pi^2 + 1 \not \in \one spec^I(\nilmfld1,g).$$
\medskip

However 
for $\tau = \pm\frac1{2}\beta_1,$
$det(E_\tau - (\pi^2+1)I_7)=0.$ Thus $\pi^2+1$ is an eigenvalue for the 
Laplacian acting on $\H_{\pm\frac1{2}\beta_1} \otimes \Lambda^1(\mm^*).$
Indeed, the eigenspace of $\pi^2+1$ in $\H^I\otimes\Lambda^1(\mm^*)$ is $$span_\Cplx\{F_{\frac{1}{2}\beta_1}\otimes\zeta_2\ , \ F_{-\frac{1}{2}\beta_1}\otimes\zeta_2 \},$$ 
which has dimension 2.  
\medskip

Thus $\pi^2 + 1 \in \one spec^I(\nilmfld2,g)$ with multiplicity 2, as desired.
\enddemo
\bigskip

The proof of Proposition 4.7 is now complete.
\endpf
\enddemo
\bigskip

\bigskip


\proclaim{Proposition 4.9}
The nilmanifolds $(\nilmfld1, g)$ and $(\nilmfld2,g)$ as presented in Example IV 
are not isospectral on one-forms.  
\endproclaim
\bigskip

\pagebreak

\demo{Proof of Proposition 4.9}
\medskip

Again using the notation of Section 2,
for $i=1,2,$ let $\Tau_{i}$  be a subset of $\mm^*$ such that $$\fnsp{i} \cong {\underset{\tau\in\Tau_{i}}\to\bigoplus} {m_i(\tau)}\H_\tau.$$
\medskip 

Let $\Tau_i = \Tau'_i \cup \Tau''_i \cup \Tau'''_i,$
where $$\align \Tau'_{i} &= \{\tau \in \Tau_{i} :  \tau(\m1) \equiv 0\}, \\
\Tau''_{i} &= \{\tau \in \Tau_{i} :  \tau(\m2) \equiv 0, \tau(\m1) \not \equiv 0\},\\
\Tau'''_{i} &= \{\tau \in \Tau_{i} :  \tau(\m2) \not \equiv 0\}.
\endalign$$
As in the proof of Proposition 4.7, decompose the representation spaces and 
spectrum accordingly.
\medskip 
 
By Lemma 3.4, the representations of $G$ on $\H'''_1$ and $\H'''_2$ are unitarily 
equivalent and 
$$\one \spec'''(\nilmfld1, g) =  \one \spec'''(\nilmfld2, g).$$
\medskip

A calculation almost identical to Step 2 in the proof of Proposition 4.7 shows 
that the representations of $G$ on $\H''_1$ and $\H''_2$ are unitarily equivalent.
Thus $$\one \spec''(\nilmfld1, g) = \one \spec''(\nilmfld2,g).$$
\medskip

It remains to show that $\one \spec'(\nilmfld1,g) \neq \one \spec'(\nilmfld2,g).$
The proof of this corresponds to Step 4 in the proof of Proposition 4.7.
\medskip

Let $\{\alpha_1, \beta_1, \beta_2, \gamma , \omega \}$ be the dual to the 
orthonormal basis
$\{X_1,Y_1, Y_2, Z, W\}$ given in Example IV. 
For $\tau \in \Tau'_i, \  \tau = A_1\alpha_1 +  B_1\beta_1 + B_2\beta_2$  for 
some $A_1, B_1, 
B_2 \in \R.$   We again use Proposition A.2. 
\medskip

Now $\mbar_\tau = \mbar.$  Hence 
$\tau(log\GGGbar{i} \cap \mbar_\tau) \subset \Z$ if and only if 
$\tau(log \GGGbar{i}) \subset \Z.$  
\medskip

Thus, $\tau \in \Tau^I_1$ if and only if 
$$A_1 \in \frac1{2}\Z \quad \text{ and } \quad B_1,B_2 \in \Z,\tag{$*$}$$  and,
$\tau \in \Tau^I_2$ if and only if $$A_1,B_2 \in \Z \quad  \text{ and } \quad B_1 \in \frac1{2}\Z.\tag{$**$}$$ 
\medskip

Let $\H_\tau$ be the associated representation space of $\pi_\tau.$ 
 As in the proof of Proposition 4.7, 
$\H_\tau \otimes \Lambda^1(\mm^*) = \F \otimes \Cplx \Lambda^1(\mm^*)$ where 
$$\split \F(exp(x_1X_1)exp(y_1&Y_1)exp(y_2Y_2)exp(zZ)exp(wW))\\
& = exp\{2\pi i \tau(x_1X_1 + y_1Y_1 + y_2Y_2) \}.\endsplit$$
\medskip

Let $\F \otimes \mu \in \F \otimes \Cplx \Lambda^1(\mm^*)$
with $\mu = a_1\alpha_1 + b_1\beta_1 + b_2\beta_2 + z\zeta + w\omega$ 
for some $a_1, b_1, b_2, z, w \in \Cplx.$
\medskip

Then 
$\Delta \F = 4\pi^2S^2\F$
where $S^2 = A_1^2 + B_1^2 + B_2^2.$ And 
$\Delta\alpha_1 = \Delta\beta_1 = \Delta\beta_2 = 0,$
$\Delta\zeta = \zeta$ and $\Delta\omega=2\omega.$ 
We also obtain the following chart:
\medskip

{\eightpoint
$$
\vbox{\offinterlineskip
\halign{\strut\vrule#&\quad#\hfil\quad&&\vrule#&\quad\hfil#\hfil\quad\cr
\noalign{\hrule}
&&&&&&&&&&&&\cr
&$\nabla_U\mu$&&$\alpha_1$&&$\beta_1$&&$\beta_2$&&$\zeta$&&$\omega$&\cr
&&&&&&&&&&&&\cr
\noalign{\hrule}
\noalign{\hrule}
&&&&&&&&&&&&\cr
&$X_1$&&0&&$\frac1{2}\zeta$&&0&&$-\frac1{2}\beta_1+\frac1{2}\omega$&&$-\frac1{2}\zeta$&\cr
&&&&&&&&&&&&\cr
\noalign{\hrule}
&&&&&&&&&&&&\cr
&$Y_1$&&$-\frac1{2}\zeta$&&0&&$\frac1{2}\omega$&&$\frac1{2}\alpha_1$&&$-\frac1{2}\beta_2$&\cr
&&&&&&&&&&&&\cr
\noalign{\hrule}
&&&&&&&&&&&&\cr
&$Y_2$&&0&&$-\frac1{2}\omega$&&0&&0&&$\frac1{2}\beta_1$&\cr
&&&&&&&&&&&&\cr
\noalign{\hrule}
}}
$$
}
\medskip 

Using (A.1) and the above information,   
if we let 
$$E_\tau=\left(\matrix 
4\pi^2S^2 & 0 & 0 & -2\pi i B_1 & 0\\
0 & 4 \pi^2S^2 & 0 & 2\pi i A_1 & -2 \pi i B_2\\
0 & 0 & 4 \pi^2S^2 & 0 & 2\pi i B_1\\
2\pi i B_1 & -2 \pi i A_1 & 0 & 4 \pi^2 S^2 + 1 & 2 \pi i A_1 \\
0 & 2 \pi i B_2 & -2 \pi i B_1 & -2 \pi i A_1 & 4 \pi^2S^2 + 2 
\endmatrix
\right),$$
then  $\Delta(\F \otimes \mu) = \lambda(\F \otimes \mu)$ if and only if $\lambda$ 
is an eigenvalue of the matrix  $E_\tau.$
\medskip

We calculate necessary conditions on $\tau=A_1\alpha_1+B_1\beta_1+B_2\beta_2$ so 
that $\pi^2+1$ is an eigenvalue of $E_\tau.$  As we assumed $\tau\in\Tau'_1$ or 
$\tau\in\Tau'_2,$  we 
know  $A_1, B_1, B_2 \in \Q.$
If $det(E_\tau - (\pi^2+1)I_5) = 0,$ then $\pi$ is the root of a polynomial with 
rational coefficients.  However $\pi$ is transcendental.  Thus the coefficients 
of the powers of $\pi$ must be zero.  
\medskip

A straightforward calculation shows that $\pi^{10}$ is the highest power 
of $\pi$ occurring in the polynomial, and the coefficient of $\pi^{10}$ 
is equal to $(4S^2-1)^5.$  Thus if $\pi^2+1$ is an eigenvalue of 
$E_\tau,$ then $S^2 = \frac1{4}.$
Recall $S^2 = A_1^2 + B_1^2 + B_2^2.$
\medskip

For $\tau$ in $\Tau'_1,$ \  $S^2=\frac1{4}$ if and only 
if (see ($*$)) $\tau = \pm\frac1{2}\alpha_1.$
And for $\tau$ in $\Tau'_2,$ \  $S^2=\frac1{4}$ if and only 
if (see ($**$)) $\tau = \pm\frac1{2} \beta_1.$
\medskip

For $\tau = \pm\frac1{2}\alpha_1,$
$det(E_\tau - (\pi^2+1)I_5) = 0.$  The eigenspace of 
$\pi^2+1$ in $\H'\otimes\Lambda^1(\mm^*)$ is
$$span_\Cplx \{F_{\frac1{2}\alpha_1}\otimes(\pi i\beta_1 + \zeta + \pi i\omega), \  
F_{-\frac{1}{2}\alpha_1}\otimes(\pi i\beta_1-\zeta+\pi i\omega)\},$$
which has dimension 2.
Thus $\pi^2 + 1 \in \one spec'(\nilmfld1, g)$ with multiplicity two.
\medskip

However, for $\tau = \pm\frac1{2}\beta_1,$ \quad 
$det(E_\tau - (\pi^2+1)I_5) \neq 0.$  
Thus 
$\pi^2 + 1$ cannot occur in $\one spec'(\nilmfld2, g)$  
and $$\one spec'(\nilmfld1, g) \neq \one spec'(\nilmfld2, g),$$
as desired.
\medskip

The proof of Proposition 4.9 is now complete. 
$\endpf$
\enddemo 
\bigskip

\bigskip


\proclaim{Proposition 4.11}
The nilmanifolds $(\nilmfld1, g)$ and $(\nilmfld2,g)$ as presented in 
Example V are not 
isospectral on one-forms.  
\endproclaim
\bigskip

\demo{Proof of Proposition 4.11}
\bigskip

Using the notation of Section 2, for 
$i=1,2,$ let $\Tau_{i}$ be a subset of $\mm^*$ such that 
$$\fnsp{i} \cong {\underset{\tau\in\Tau_{i}}\to\bigoplus} {m_i(\tau)}\H_\tau.$$
\medskip 

Let $\Tau_i = \Tau^{I}_i \cup \Tau^{II}_i \cup \Tau^{III}_i \cup \Tau^{IV}_i$
be defined as in the proof of Proposition 4.7,
and decompose the representation spaces and spectrum accordingly.
\medskip 

By Lemma 3.4, the representations of $G$ on $\H^{IV}_1$ and $\H^{IV}_2$ are 
unitarily equivalent
and  $$\one \spec^{IV}(\nilmfld1, g) =  \one \spec^{IV}(\nilmfld2, g).$$
\medskip

A calculation almost identical to that in Step 2 in the proof of 
Proposition 4.7 shows that the representations of $G$ on $\H^{III}_1$ 
and $\H^{III}_2$ are unitarily equivalent. 
Thus $$\one \spec^{III}(\nilmfld1, g) = \one \spec^{III}(\nilmfld2,g).$$
\medskip

A calculation almost identical to that in Step 3 in the proof of 
Proposition 4.7 shows that the representations of $G$ on $\H^{II}_1$ 
and $\H^{II}_2$ are unitarily equivalent. 
Thus $$\one \spec^{II}(\nilmfld1, g) = \one \spec^{II}(\nilmfld2,g).$$
\medskip

It remains to show that  $\one spec^I(\nilmfld1, g) \neq \one spec^I(\nilmfld1, g).$
\medskip

For $\tau \in \Tau^I_1$ or $\tau \in \Tau^I_2,\  \tau(\g1)\equiv 0. $  We again use 
Proposition A.2 to calculate the irreducible representations occurring here.  Let 
$\{\epsilon_1, \cdots, \epsilon_7\}$ be the dual to the orthonormal basis 
$\{E_1, \cdots, E_7\}$ given in Example V. Then
$\tau = A_1\epsilon_1 + A_2\epsilon_2 +  A_3\epsilon_3 + A_4\epsilon_4$ for some 
$A_1, \cdots, A_4 \in \R.$ 
\medskip

Now $\mbar_\tau = \mbar.$  Hence 
$\tau(log\GGGbar{i} \cap \mbar_\tau) \subset \Z$ if and only if 
$\tau(log \GGGbar{i}) \subset \Z.$  
\medskip

Thus, $\tau \in \Tau^I_1$ if and only if 
$$2A_1+A_2-\frac{1}{4}A_3+\frac{1}{2}A_4\in\Z, \quad 2A_2+\frac{1}{2}A_3\in\Z,\quad\text{and}\quad A_3,A_4 \in \Z.\tag{$*$}$$  And
$\tau \in \Tau^I_2$ if and only if  $$2A_1+A_2+\frac{1}{4}A_3+\frac{1}{2}A_4\in\Z,\quad  2A_2+\frac{1}{2}A_3\in\Z,\quad\text{and}\quad  A_3,A_4 \in \Z. \tag{$**$}$$
Note that the only distinction between the two conditions is in the sign of $\frac{1}{4}A_3.$ 
\medskip

Let $\H_\tau$ be the associated representation space of $\pi_\tau.$  
As in the proof 
of Proposition 4.7, 
$\H_\tau \otimes \Lambda^1(\mm^*) = \F \otimes \Cplx \Lambda^1(\mm^*)$ 
where
$$\split \F(exp(e_1E_1)&exp(e_2E_2)exp(e_3E_3)exp(e_4E_4)exp(e_5E_5)exp(e_6E_6)exp(e_7E_7)) \\
&= exp\{2\pi i \tau(e_1E_1 + e_2E_2 + e_3E_3 + e_4E_4) \}.\endsplit$$
\medskip

Let $\F \otimes \mu \in \F \otimes \Cplx \Lambda^1(\mm^*)$
with  $\mu = a_1\epsilon_1 + a_2\epsilon_2 + a_3\epsilon_3 
+ a_4\epsilon_4 + a_5\epsilon_5 + a_6\epsilon_6 + a_7\epsilon_7$ 
for some $a_i \in \Cplx,$ for $i=1, \dots, 7.$
As in the proof of Proposition 4.7,
$\Delta \F = 4\pi^2S^2\F$
where 
$S^2 = A_1^2 + A_2^2 + A_3^2 + A_4^2.$
Also, $\Delta\epsilon_1 = \Delta\epsilon_2 = \Delta\epsilon_3 = \Delta\epsilon_4 = 0,$
$\Delta\epsilon_5 = 2\epsilon_5,$ $\Delta\epsilon_6=\epsilon_6+\frac{1}{4}\epsilon_7,$  and 
$\Delta\epsilon_7 = \frac{1}{4}\epsilon_6+(3+\frac{1}{256}+\frac{3}{16})\epsilon_7.$ 
And we obtain the following chart:
\bigskip

{\eightpoint
\centerline{
\vbox{\offinterlineskip
\halign{\strut\vrule#& #\hfil &&\vrule#& \hfil#\hfil \cr
\noalign{\hrule}
&&&&&&&&&&&&&&&&\cr
&\ $\nabla_U\mu\ $&&$\epsilon_1$&&$\epsilon_2$&&$\epsilon_3$&
&$\epsilon_4$&&$\epsilon_5$&&$\epsilon_6$&
&$\epsilon_7$&\cr
&&&&&&&&&&&&&&&&\cr
\noalign{\hrule}
\noalign{\hrule}
&&&&&&&&&&&&&&&&\cr
&\ \  $E_1$&&0&&$-\frac{1}{32}\epsilon_7$&&$\frac1{2}\epsilon_5+\frac{1}{8}\epsilon_7$&
&$\frac1{2}\epsilon_6+\frac{1}{8}\epsilon_7$&&$-\frac1{2}\epsilon_3+\frac1{2}\epsilon_7$&
&$-\frac1{2}\epsilon_4$&
&$\frac{1}{32}\epsilon_2-\frac{1}{8}\epsilon_3-\frac{1}{8}\epsilon_4-\frac{1}{2}\epsilon_5$&\cr
&&&&&&&&&&&&&&&&\cr
\noalign{\hrule}
&&&&&&&&&&&&&&&&\cr
&\ \  $E_2$&&$\frac{1}{32}\epsilon_7$&&0&&0&&$\frac1{2}\epsilon_5-\frac{1}{8}\epsilon_7$&
&$-\frac{1}{2}\epsilon_4$&&$\frac{1}{2}\epsilon_7$&
&$-\frac{1}{32}\epsilon_1+\frac{1}{8}\epsilon_4-\frac1{2}\epsilon_6$&\cr
&&&&&&&&&&&&&&&&\cr
\noalign{\hrule}
&&&&&&&&&&&&&&&&\cr
&\ \  $E_3$&&$-\frac1{2}\epsilon_5-\frac{1}{8}\epsilon_7$&&0&&0&&$\frac1{2}\epsilon_7$&
&$\frac1{2}\epsilon_1$&&$0$&&$\frac{1}{8}\epsilon_1-\frac1{2}\epsilon_4$&\cr
&&&&&&&&&&&&&&&&\cr
\noalign{\hrule}
&&&&&&&&&&&&&&&&\cr
&\ \  $E_4$&&$-\frac1{2}\epsilon_6-\frac{1}{8}\epsilon_7$&
&$-\frac1{2}\epsilon_5+\frac{1}{8}\epsilon_7$&&$-\frac1{2}\epsilon_7$&&0&
&$\frac1{2}\epsilon_2$&&$\frac1{2}\epsilon_1$&
&$\frac{1}{8}\epsilon_1-\frac{1}{8}\epsilon_2+\frac1{2}\epsilon_3$&\cr
&&&&&&&&&&&&&&&&\cr
\noalign{\hrule}
}}
}
}
\bigskip

Using (A.1) and the above information, a straightforward calculation shows that
if we let $E_\tau $ be the skew-Hermitian matrix defined by
{\eightpoint
$$\left(\matrix 
4\pi^2S^2 & 0 & 0 & 0 &  2\pi i A_3 & 2\pi i A_4 & \pi i (-\frac{1}{8}A_2+\frac{1}{2}A_3+\frac{1}{2}A_4)\\
 & 4 \pi^2S^2 & 0 & 0 & 2 \pi i A_4 & 0 & \pi i (\frac{1}{8}A_1-\frac{1}{2}A_4)\\
 &  & 4 \pi^2S^2 & 0 & -2 \pi i A_1 & 0   & \pi i ( -\frac{1}{2}A_1 + 2 A_4)\\
 &  &  & 4 \pi^2S^2 & -2 \pi i A_2& -2 \pi i A_1 & \pi i (-\frac{1}{2}A_1 +\frac{1}{2}A_2 - 2 A_3)\\
  & &  &  &  4 \pi ^2 S^2 + 2 &0 & -2 \pi i A_1 \\
 &  &  &  &  &  4 \pi ^2 S^2 + 1 & - 2 \pi i A_2 +\frac{1}{4}\\
 & & & & & & 4 \pi^2S^2 + \frac{817}{256}
\endmatrix
\right)$$
}
then $\Delta(\F \otimes \mu) = \lambda(\F \otimes \mu)$ if and 
only if $\lambda$ is an eigenvalue of the matrix  $E_\tau.$
\medskip

Let $$\lambda = \frac{17}{4}\pi^2+1+\sqrt{\frac{17}{4}\pi^2+1}.$$
\medskip

We now calculate necessary conditions on $\tau = A_1\epsilon_1 + 
A_2\epsilon_2 +  A_3\epsilon_3 + A_4\epsilon_4$
so that $\lambda$ is an eigenvalue of $E_\tau.$
As $\tau\in\Tau^I_1$ or $\tau\in\Tau^I_2,$ we know  $A_1, A_2, A_3, A_4 \in \Q.$
\medskip

By a computation using Maple or Mathematica, if $det(E_\tau - \lambda I_7) = 0,$
then $\pi$ is the root of a polynomial with rational coefficients. 
However $\pi$ is transcendental.  Thus the coefficients of the powers of $\pi$ must be zero.  
\medskip

A straightforward calculation shows that $\pi^{14}$ is the highest power of 
$\pi$ occurring in the polynomial, and the coefficient of $\pi^{14}$ 
is equal to $(4S^2-\frac{17}{4})^7.$  Thus if $\lambda$ is an 
eigenvalue of $E_\tau,$ then $S^2 = \frac{17}{16}.$  Recall that 
$S^2=A_1^2+A_2^2+A_3^2+A_4^2.$
\medskip

For $\tau$ in $\Tau^I_1,$ \   $S^2 = \frac{17}{16}$ if and only 
if (see ($*$))  $\tau = \pm (\frac{1}{4}\epsilon_2+\epsilon_3)$ 
or $\tau = \pm \frac{1}{4}\epsilon_1 \pm \epsilon_4.$
And for $\tau$ in $\Tau^I_2,$ \   $S^2 = \frac{17}{16}$ if and 
only 
if (see ($**$)) $\tau = \pm (\frac{1}{4}\epsilon_2-\epsilon_3)$ 
or $\tau = \pm \frac{1}{4}\epsilon_1 \pm \epsilon_4.$  
Note that the only difference is in the sign of $\epsilon_3.$
\medskip

For $\tau = \pm (\frac{1}{4}\epsilon_2+\epsilon_3)$ or 
$\tau = \pm \frac{1}{4}\epsilon_1 \pm \epsilon_4,$ 
a calculation using Maple or Mathematica shows that 
$det(E_\tau - \lambda I_7) \neq 0.$  
Thus $\frac{17}{4}\pi^2+1+\sqrt{\frac{17}{4}\pi^2+1}\not \in \one spec^I(\nilmfld1,g).$
\medskip

However, 
for $\tau =  \pm (\frac{1}{4}\epsilon_2-\epsilon_3),$
$det(E_\tau - \lambda I_7)=0.$ Thus  
$\frac{17}{4}\pi^2+1+\sqrt{\frac{17}{4}\pi^2+1}$ 
is an eigenvalue for the Laplacian acting on 
$\H_{\pm (\frac{1}{4}\epsilon_2-\epsilon_3)} \otimes \Lambda^1(\mm^*),$ and 
$\frac{17}{4}\pi^2+1+\sqrt{\frac{17}{4}\pi^2+1} \in \one spec^I(\nilmfld2,g).$
\medskip  
 
Consequently $$1\hy spec^I(\nilmfld1, g) \neq 1\hy spec^I(\nilmfld2, g),$$
as desired. 
\medskip

The proof of Proposition 4.11 is now complete.
$\endpf$
\enddemo

\bigskip

\bigskip


\enddocument